\documentclass{article}
\usepackage{latexsym,amssymb,amsmath}
\newtheorem{theorem}{Theorem}[section]
\newtheorem{lemma}[theorem]{Lemma}

\newtheorem{definition}[theorem]{Definition}
\newtheorem{proposition}[theorem]{Proposition}
\newtheorem{corollary}[theorem]{Corollary}
\newtheorem{remark}[theorem]{Remark}
\newtheorem{dfn}[theorem]{Definition}

\newtheorem{cor}[theorem]{Corollary}
\newtheorem{thm}[theorem]{Theorem}

\newtheorem{thmo}[theorem]{Theorem}

\newtheorem{deff}[theorem]{Definition}

\newtheorem{prop}[theorem]{Proposition}

\newcommand{\nel}{\mathbb{N}  }

\newcommand{\TP}{OT(\psi)}
\newcommand{\TT}{\mathrm{OT}_{\!_{\mathfrak X}}(\vartheta)}

\newcommand{\Pie}{\Pi^{1}_{2}}

\newcommand{\Ts}[3]{\provx{T^*_{\!_Q}}{#1}{#2}{#3}}
\newcommand{\Tsr}[2]{\provx{T^*_{\!_Q}}{#1}{\varrho}{#2}}

\newcommand{\Tst}{T^*_{\!_Q}}

\newcommand{\Om}{\Omega}

\newcommand{\kOm}{<_{\Omega}}

\newcommand{\om}{\omega}
\newcommand{\al}{\alpha}
\newcommand{\be}{\beta}
\newcommand{\ga}{\gamma}
\newcommand{\de}{\delta}

\newcommand{\vt}{\vartheta} 
\newcommand{\vtt}{\vartheta_{\!_{\!\mathfrak X}}}
\newcommand{\mq}{\overline{m}}

\newcommand{\Gt}{G_{\vt }}

\newcommand{\und}{\;\wedge\;}
\newcommand{\AO}{D_{\Omega}}

\newcommand{\proof}{{\bf Proof}: }
\newcommand{\dX}{|{\mathfrak X}|}

\newcommand{\EE}{{\mathfrak E}}
\newcommand{\CC}{C_{_{\!\mathfrak X}}}
\newcommand{\kvt}{<}
\newcommand{\kgvt}{\leq}
\newcommand{\vd}{\vert \vtt \vert}
\newcommand{\kgrvt}{\geq}
\newcommand{\ACA}{{\mathbf{ACA}}}
\newcommand{\BI}{{\mathbf{BI}}}
\newcommand{\Acc}{{\mathrm{Acc}}}
\newcommand{\Mm}{{\mathrm{M}}}
\newcommand{\RCA}{{\mathbf{RCA}}}
\newcommand{\ATR}{{\mathbf{ATR}}}
\newcommand{\WO}{{\mathrm{WO}}}
\newcommand{\Prog}{{\mathrm{Prog}}}
\newcommand{\Mod}{{\mathfrak{A}}}
\newcommand{\WF}{{\mathrm{WF}}}
\newcommand{\TI}{{\mathrm{TI}}}
\newcommand{\prf}{{\bf Proof}. }
\newcommand{\Ii}{{\mathcal I}}
\newcommand{\fXO}{{\mathfrak X}}
 \newcommand{\bbM}[1]{{\mathbb #1}}
 \newcommand{\DQ}{\bar{Q}}
 \newcommand{\MQ}{{\mathcal D}_Q}
 \newcommand{\Sr}{{\mathbf S}}
 \newcommand{\SAC}{\Sigma^1_1\mbox{-}{\mathbf{AC}}}
 \newcommand{\calD}{{\mathcal D}}
 \newcommand{\astr}{*}
\newcommand{\gdw}{\leftrightarrow}
\newcommand{\WOP}{{\mathbf{WOP}}}
\newcommand{\fX}{{\mathfrak X}}
\newcommand{\paar}[1]{\langle #1\rangle}

 \def\provx#1#2#3#4{
\setbox1=\hbox{\kern1.5pt$\scriptstyle#3$}
\def\zeichen{#2}
\ifx\zeichen\empty\setbox0=\hbox to .75em{}\else\setbox0=\hbox
{\kern1.5pt$\scriptstyle#2$}\fi
\dimen1=\dp0 \ifdim \dimen1=0pt
\advance \dimen1 by 1.5ex \else \advance \dimen1 by 1.2ex
\fi\dimen3=2ex\dimen4=.5ex\ifdim \wd0<\wd1 \dimen2=\wd1 \else \dimen2=\wd0
\fi\hbox{$#1\hskip 5pt minus5pt\vrule height\dimen3
depth\dimen4\raise\dimen1\copy0\hskip-1\wd0 \lower\ht1
\copy1\hskip-1\wd1\vrule width\dimen2 height.7ex depth-.6ex\hskip3pt
minus1.5pt#4\hskip2pt plus2pt minus2pt$}}

\def\prov#1#2#3{
\setbox1=\hbox{\kern1.5pt$\scriptstyle#2$}
\def\zeichen{#1}
\ifx\zeichen\empty\setbox0=\hbox to .75em{}\else\setbox0=\hbox
{\kern1.5pt$\scriptstyle#1$}\fi
\dimen1=\dp0
\ifdim \dimen1=0pt
\advance \dimen1 by 1.5ex \else \advance \dimen1 by 1.2ex
\fi\dimen3=2ex\dimen4=.5ex\ifdim \wd0<\wd1 \dimen2=\wd1 \else \dimen2=\wd0
\fi\hbox{\hskip0pt plus 4pt
$\vrule height\dimen3
depth\dimen4\raise\dimen1\copy0\hskip-1\wd0
\lower\ht1\copy1\hskip-1\wd1\vrule width\dimen2 height.7ex depth-.6ex
\hskip3pt minus1.5pt#3\hskip2pt plus2pt minus2pt$}}

\def\prv#1#2{
\setbox1=\hbox{\kern1.5pt$\scriptstyle#2$}
\ifx\zeichen\empty\setbox0=\hbox to .75em{}\else\setbox0=\hbox
{\kern1.5pt$\scriptstyle#1$}\fi
\dimen1=\dp0 \ifdim \dimen1=0pt
\advance \dimen1 by 1.5ex \else \advance \dimen1 by 1.2ex
\fi\dimen3=2ex\dimen4=.5ex\ifdim \wd0<\wd1 \dimen2=\wd1 \else \dimen2=\wd0
\fi\hbox{\hskip.5em$\vrule height\dimen3
depth\dimen4\raise\dimen1\copy0\hskip-1\wd0
\lower\ht1\copy1\hskip-1\wd1\vrule width\dimen2 height.7ex depth-.6ex
\hskip3pt minus1.5pt$}}

\mathchardef\str='1066
\def\negprov#1#2#3{
\setbox1=\hbox{\kern1.5pt$\scriptstyle#2$}
\setbox4=\hbox{$\str$}
\def\zeichen{#1}
\ifx\zeichen\empty\setbox0=\hbox to 1em{}\else\setbox0=\hbox
{\kern1.5pt$\scriptstyle#1$}\fi
\dimen1=\dp0
\ifdim \dimen1=0pt
\advance \dimen1 by 1.5ex \else \advance \dimen1 by 1.2ex
\fi\dimen3=2ex\dimen4=.5ex\ifdim \wd0<\wd1 \dimen2=\wd1 \else \dimen2=\wd0
\fi\hbox{\hskip.5em$\kern-1.9pt\raise1pt\copy4\kern-\wd4\kern1.9pt\vrule height\dimen3
depth\dimen4\raise\dimen1\copy0\hskip-1\wd0
\lower\ht1\copy1\hskip-1\wd1\vrule width\dimen2 height.7ex depth-.6ex
\hskip3pt minus1.5pt#3\hskip2pt plus2pt minus2pt$}}

\def\goed#1{\setbox5=\hbox{$#1$}\dimen1=.25em \dimen2=\dimen1 \advance \dimen2
by -1pt\hbox{\raise.65\ht5 \hbox{\vrule height.5\ht5 depth0pt width.4pt\vrule
height.5\ht5 width\dimen1 depth-.48\ht5}\kern-\dimen2\copy5\kern-\dimen2
\raise.65\ht5 \hbox{\vrule height .5\ht5 width\dimen1 depth-.48\ht5\vrule
height.5\ht5 depth 0pt width.4pt}\hskip4pt plus2pt minus2pt}}

\def\mod#1#2{
\def\zeichen{#1}
\hbox{\hskip 2pt plus3pt minus 2pt\vrule width.5pt height2ex depth.5ex
\vbox{\ifx\zeichen\empty\hbox to .75em{}\else
\hbox{\kern1.5pt $\scriptstyle#1$}\fi
\kern2pt
\hrule
\kern1.7pt
\hrule\kern1.7pt}
\hskip3pt minus 2pt$#2$}\hskip2pt
plus3pt minus2pt}

\def\notmod#1#2{\hbox{\hskip 2pt plus 3pt minus 3pt\vrule width.5pt
height2ex depth.5ex
\vbox{\hbox{\kern1.5pt $\scriptstyle#1$}\kern3pt
\setbox0=\hbox{\kern2pt$\scriptstyle/$}
\hrule
\kern-1.7pt
\copy0
\kern-\ht0
\kern 1.7pt
\hrule\kern1.7pt}n
\hskip3pt minus 2pt$#2$}\hskip2pt
plus3pt minus2pt}
\def\sq{\hbox{\rlap{$\sqcap$}$\sqcup$}}
\def\qed{\ifmmode\sq\else{\unskip\nobreak\hfil\penalty50\hskip1em\null
\nobreak\hfil\sq\parfillskip=0pt\finalhyphendemerits=0\endgraf}\fi\medskip}

\def\lleq{\hbox{\hskip3pt minus3pt\kern1pt\lower4pt
\vbox{\hbox{$\scriptstyle\ll$}
\kern-7pt\hbox{\kern1pt$\scriptstyle=$}}\hskip3pt minus 3pt}}

\mathchardef\res='1152
\mathchardef\qin='1062
\mathchardef\qprec='1036
\mathchardef\qless='474
\mathchardef\dpkt='72

\def\A{{\rm\kern.22em
   \vrule width.02em
       height0.5ex depth 0ex
   \kern-.24em A}}
\def\B{{\rm I\kern-.25em B}}
\def\C{{\rm\kern.24em
   \vrule width.02em
       height1.4ex depth-.05ex
   \kern-.26em C}}
\def\D{{\rm I\kern-.25em D}}
\def\E{{\rm I\kern-.25em E}}
\def\F{{\rm I\kern-.25em F}}
\def\G{{\rm\kern.24em
   \vrule width.02em
       height1.4ex depth-.05ex
   \kern-.26em G}}
\def\H{{\rm I\kern-.25em H}}
\def\I{{\rm I\kern-.25em I}}
\def\J{{\rm\kern.19em
   \vrule width.02em
       height1.47ex depth 0ex
   \kern-.21em J}}
\def\K{{\rm I\kern-.25em K}}
\def\L{{\rm I\kern-.25em L}}
\def\M{{\rm I\kern-.23em M}}
\def\N{{\rm I\kern-.23em N}}
\def\O{{\rm\kern.24em
   \vrule width.02em
       height1.4ex depth-.05ex
   \kern-.26em O}}
\def\P{{\rm I\kern-.25em P}}
\def\Q{{\rm\kern.24em
   \vrule width.02em
       height1.4ex depth-.05ex
   \kern-.26em Q}}
\def\R{{\rm I\kern-.25em R}}
\def\S{{\rm\kern.18em
   \vrule width.02em
       height1.47ex depth-.9ex
\kern.12em
\vrule width.02em
    height0.7ex depth 0ex
\kern-.34em S}}
\def\T{{\rm\kern.45em
   \vrule width.02em
       height1.47ex depth 0ex
   \kern-.47em T}}
\def\U{{\rm\kern.30em
   \vrule width.02em
       height1.47ex depth-.05ex
   \kern-.32em U}}
\def\V{{\rm\kern.27em
   \vrule width.02em
       height1.47ex depth-.8ex
   \kern-.29em V}}
\def\W{{\rm\kern.25em
   \vrule width.02em
       height1.47ex depth-.9ex
   \kern.34em
   \vrule width.02em
       height1.47ex depth-.9ex
   \kern-.63em W}}
\def\X{{\rm\kern.30em
   \vrule width.02em
       height1.47ex depth-1ex
   \kern.12em
   \vrule width.02em
       height0.4ex depth 0ex
   \kern-.46em X}}
\def\Y{{\rm\kern.25em
   \vrule width.02em
       height1.0ex depth 0ex
   \kern-.27em Y}}
\def\Z{{\rm\kern.26em
   \vrule width.02em
       height0.5ex depth 0ex
   \kern.04em
   \vrule width.02em
       height1.47ex depth-1ex
   \kern-.34em Z}}

   \newcommand{\On}{{\mathrm{ON}}}
   \newcommand{\AP}{{\mathrm{AP}}}
   \newcommand{\Ep}{{\mathrm{E}}}

\begin{document}
\title{Well ordering principles and bar induction}
\author{Michael Rathjen and Pedro Francisco Valencia Vizca\'ino\\
 {\small {\it School of Mathematics,
  University of Leeds}} \\
{\small {\it Leeds, LS2 JT, England}}}
\date{}
\maketitle

\begin{abstract}In this paper we show that the existence of $\omega$-models of bar induction is equivalent to
the principle saying that applying the Howard-Bachmann operation to any well-ordering
yields again a well-ordering.
\\[1ex]
MSC 03B30 03F05 03F15 03F35 03F35
\end{abstract}

\section{Introduction}
This paper  will be concerned with a particular $\Pi^1_2$ statement
of the form
 \begin{eqnarray}\label{1}\WOP(f):&&\;\,\,\;\;\;\forall X\,[\WO(\fX)\rightarrow \WO(f(\fX))]\end{eqnarray} where $f$ is a standard proof-theoretic function from ordinals to ordinals and $\WO(\fX)$ stands for `$\fX$ is a well-ordering'.
 There are by now several examples of functions $f$ familiar from proof theory where the statement
 $\WOP(f)$ has turned out to be equivalent to one of the theories of reverse
 mathematics over a weak base theory (usually ${\mathbf{RCA}}_0$). The first explicit example appears to be due to Girard
 \cite[5.4.1 theorem] {girard87} (see also \cite{Hirst}). However, it is also implicit in Sch\"utte's
 proof of cut elimination for $\omega$-logic \cite{sch51} and ultimately has its roots in Gentzen's work, namely in his first unpublished consistency proof\footnote{The original German version was finally published in 1974 \cite{gentzen74}. An earlier English
translation appeared in 1969 \cite{gentzen69}.},  where he introduced the notion of a ``Reduziervorschrift" \cite[p. 102]{gentzen74} for a sequent. The latter is
a well-founded tree built bottom-up via ``Reduktionsschritte", starting with the given sequent
and passing up from conclusions to premises until an axiom is reached.

 \begin{thmo}\label{girard}
Over ${\mathbf{RCA}}_0$ the following are equivalent:
 \begin{itemize}
  \item[(i)] Arithmetical comprehension
  \item[(ii)] $\forall {\fX}\;[\WO(\fX)\rightarrow \WO(2^{\fX})]$.
  \end{itemize}
  \end{thmo}
  Another characterization from \cite{girard87}, Theorem 6.4.1, shows that arithmetical comprehension is equivalent to Gentzen's Hauptsatz (cut elimination)  for $\omega$-logic.
  Connecting statements of form (\ref{1}) to cut elimination theorems for
  infinitary logics will also be a major tool in this paper.

 There are several more recent examples of such equivalences that have been proved by
 recursion-theoretic as well proof-theoretic methods.
 These results give characterizations of the form (\ref{1}) for the theories $\mathbf{ACA}_0^+$ and $\mathbf{ATR}_0$, respectively, in terms of familiar proof-theoretic functions.
  $\mathbf{ACA}_0^+$ denotes the theory $\mathbf{ACA}_0$ augmented
  by an axiom asserting that for any set $X$ the $\omega$-th jump in $X$ exists while $\mathbf{ATR}_0$ asserts the existence of sets constructed by transfinite iterations of arithmetical comprehension. 
$\alpha\mapsto\varepsilon_{\alpha}$ denotes the usual $\varepsilon$ function while $\varphi$ stands for the two-place Veblen function
familiar from predicative proof theory (cf. \cite{sc}).
Definitions of the familiar subsystems of reverse mathematics can be found in \cite{SOSA}.

  \begin{thmo}\label{MM}  {\em (Afshari, Rathjen \cite{AR}; Marcone, Montalb\'an \cite{MM})} Over ${\mathbf{RCA}}_0$ the following are equivalent:
 \begin{itemize}
  \item[(i)] $\mathbf{ACA}_0^+$

  \item[(ii)] $\forall {\fX}\;[\WO(\fX)\rightarrow
\WO(\varepsilon_{\fX})]$.
  \end{itemize}
  \end{thmo}

   \begin{thmo}\label{FMW} {\em (Friedman \cite{FMW}; Rathjen, Weiermann \cite{rathjen-weiermann93}; Marcone, Montalb\'an \cite{MM})} Over ${\mathbf{RCA}}_0$ the following are equivalent:
 \begin{itemize}
  \item[(i)] $\mathbf{ATR}_0$

  \item[(ii)] $\forall {\fX}\;[\WO(\fX)\rightarrow
\WO(\varphi{\fX}0)]$.
  \end{itemize}
  \end{thmo}
  There is often
another way of characterizing  statements of the form (\ref{1}) by means of
the notion of countable coded $\omega$-model.

\begin{deff}\label{omegamod}{\em Let $T$ be a theory in the language of second order arithmetic, $\mathcal{L}_2$. A {\em countable coded $\omega$-model of $T$} is a set $W\subseteq {\mathbb N}$, viewed as encoding the $\mathcal{L}_2$-model
  $${\mathbb M}=({\mathbb N},{\mathcal S},\in,+,\cdot,0,1,<)$$
  with ${\mathcal S}=\{(W)_n\mid n\in{\mathbb N}\}$ such that ${\mathbb M}\models T$ when the second order quantifiers are interpreted as ranging over ${\mathcal S}$ and the first order part is interpreted in the standard way
  (where $(W)_n=\{m\mid \paar{n,m}\in W\}$ with $\paar{\,,}$ being some primitive recursive coding function).

  If $T$ has only finitely many axioms it is obvious how to express
  ${\mathbb M}\models T$ by just translating the second order quantifiers
  $QX\ldots X\ldots$ in the axioms by $Qx \ldots (W)_x\ldots$.
  If $T$ has infinitely many axioms one needs to formalize Tarski's truth definition for ${\mathbb M}$.
This definition can be made in $\mathbf{RCA}_0$ as is shown in \cite{SOSA}, Definition II.8.3 and
  Definition VII.2. Some more details will be provided in Remark \ref{tarski-truth}.

  We write $X\in W$ if $\exists n\;X=(W)_n$.
  }\end{deff}

The alternative characterizations alluded to above are as follows:

\begin{thm}\label{alle} Over ${\mathbf{RCA}}_0$ the following are equivalent:
 \begin{itemize}
  \item[(i)] $\forall {\fX}\;[\WO(\fX)\rightarrow
\WO(\varepsilon_{\fX})]$ is equivalent to the statement that every set
 is contained in a countable coded $\omega$-model of $\mathbf{ACA}$.

  \item[(ii)] $\forall {\fX}\;[\WO(\fX)\rightarrow
\WO(\varphi{\fX}0)]$ is equivalent to the statement that every set
 is contained in a countable coded $\omega$-model of
$\Delta^1_1\mbox{-}{\mathbf{CA}}$ (or $\Sigma^1_1\mbox{-}{\mathbf{DC}}$).
  \end{itemize}
  \end{thm}
  \prf See \cite[Corollary 1.8]{H-Band}. \qed

  Whereas Theorem \ref{alle} has been  established independently by recursion-theoretic  and
  proof-theoretic methods, there is also a result that has a very involved proof and so far has only been shown by
  proof theory. It connects the well-known $\Gamma$-function (cf. \cite{sc}) with the existence
  of countable coded $\omega$-models of ${\mathbf{ATR}}_0$.

   \begin{thmo}\label{R} {\em (Rathjen \cite[Theorem 1.4]{H-Band})}
    Over ${\mathbf{RCA}}_0$ the following are equivalent:
 \begin{itemize}
  \item[(i)]  $\forall {\fX}\;[\WO(\fX)\rightarrow
\WO(\Gamma_{\!\fX})]$.

  \item[(ii)] Every set
 is contained in a countable coded $\omega$-model of $\mathbf{ATR}_0$.
  \end{itemize}
  \end{thmo}
  The tools from proof theory employed in the above theorems involve search trees and
  Gentzen's cut elimination technique for infinitary logic with ordinal bounds. One could perhaps generalize and say that every cut elimination theorem in ordinal-theoretic proof theory encapsulates a theorem of this type.

  The proof-theoretic ordinal functions that figure in the foregoing theorems are all
  familiar from so-called predicative or meta-predicative  proof theory.
  Thus far a function from genuinely impredicative proof theory is missing.
  The first such function that comes to mind is of the Bachmann-Howard type.
  It was conjectured in \cite{rathjen-weiermann} (Conjecture 7.2) that the pertaining
  principle (\ref{1}) would be equivalent to the existence of countable coded $\omega$-models of bar induction, $\BI$. The conjecture is by and large true as will be shown in this paper, however,
  the relativization of the Bachmann-Howard construction allows for two different approaches, yielding principles of different strength. As it turned out, only the strongest one is equivalent
  to the existence of $\omega$-models of $\BI$.
  We now proceed to state the main result of this paper. Unexplained notions will be
  defined shortly.

  \begin{theorem}\label{main}
  Over ${\mathbf{RCA}}_0$ the following are equivalent:
 \begin{itemize}
  \item[(i)] $\mathbf{RCA}_0+\mbox{\it Every set $X$ is contained in a countable coded $\omega$-model of $\mathbf{BI}$}$.

  \item[(ii)] $\forall {\fX}\;[\WO(\fX)\rightarrow
\WO(\vartheta_{\!\fX})]$.
  \end{itemize}
  \end{theorem}
 Below we shall refer to  Theorem \ref{main} as the {\bf Main Theorem}.

\subsection{A brief outline of the paper}
Subsection \ref{s1} contains a detailed definition of the theory $\BI$.
Section \ref{s2} introduces a relativized version of the Howard-Bachmann
ordinal representation system, i.e. given a well-ordering
$\fX$, one defines a new well-ordering $\vartheta_{\!\fX}$ of Howard-Bachmann type which incorporates $\fX$. Section \ref{s3} proofs the direction
$(i)\Rightarrow(ii)$ of Theorem \ref{main}.
With section \ref{s4} the proof of  Theorem \ref{main} $(ii)\Rightarrow(ii)$ commences.
It introduces the crucial notion of a deduction chain for a given set $Q\subseteq \mathbb N$. The set of deduction chains forms a tree
${\mathcal D}_Q$. It is shown that from an infinite branch of this tree
one can construct a countable coded $\omega$-model of $\BI$ which contains $Q$. As a consequence, it remains to consider the case when
${\mathcal D}_Q$ does not contain an infinite branch, i.e.
when ${\mathcal D}_Q$ is a well-founded tree. Then the
Kleene-Brouwer  ordering of ${\mathcal D}_Q$, $\mathfrak X$, is a well-ordering and, by the
well-ordering principle (ii), $\vartheta_{\!\fX}$ is a well-ordering, too.
It will then be revealed that ${\mathcal D}_Q$ can be viewed as a skeleton of a proof
$\mathcal D^*$ of the empty sequent in an infinitary proof system $T^*_{\!_Q}$
with Buchholz' $\Omega$-rule.
However, with the help of transfinite induction over $\vartheta_{\!\fX}$
it can be shown that all cuts in $\mathcal D^*$ can be removed, yielding a
cut-free derivation of the empty sequent. As this cannot be, the final conclusion reached is
 that ${\mathcal D}_Q$ must contain an infinite branch, whence there is a
 countable coded $\omega$-model of $\BI$ containing $Q$, thereby
completing the proof of Theorem \ref{main} $(ii)\Rightarrow(i)$.

 \subsection{The theory $\BI$}\label{s1}
In this subsection we introduce the theory $\BI$.
To set the context, we fix some notations.
The language of second order arithmetic, ${\mathcal L}_2$, consists of free
numerical variables $a,b,c,d,\ldots$, bound numerical variables
$x,y,z,\ldots$, free set variables $U,V,W, \ldots,$ bound set variables
$X,Y,Z, \ldots$, the constant $0$, a symbol for each primitive recursive
function, and the symbols $=$ and $\in$ for equality in the first sort
and the elementhood relation, respectively. The {\it numerical terms} of
${\mathcal L}_2$ are built up in the usual way; $r,s,t, \ldots$ are
syntactic variables for them. {\it Formulas} are obtained from atomic
formulas $s=t$, $s \in U \mbox{ and negated atomic formulas }
\neg\,s=t, \neg\,s \in U$ by closing under $\wedge, \vee$ and
quantification $\forall x, \exists x, \forall X, \exists X$ over both sorts; so we stipulate
that formulas are in negation normal form.

The classes of $\Pie$-- and $\Sigma^{1}_{n}$--{\em formulae} are defined
as
usual (with $\Pi^{1}_{0} = \Sigma^{1}_{0} = \cup \{ \Pi^{0}_{n} : n \in
{\mathbb N} \}$). $\neg A$ is defined by de Morgan's laws; $A \to B$ stands for
$\neg A\; \vee\; B$. All theories in ${\mathcal L}_2$ will be assumed to contain
the axioms and rules of classical two sorted predicate calculus, with
equality in the first sort. In addition, it will be assumed that they
comprise the system $\ACA_{0}$. $\ACA_0$ contains all axioms of elementary
number theory, i.e. the usual axioms for $0$, $'$ (successor), the
defining equations for the primitive recursive functions, the {\it
induction axiom}
$$\forall X \,[0 \in X \;\wedge\; \forall x (x \in X \to x' \in X) \to
\forall x (x \in X)]
,$$
and all instances of {\it arithmetical comprehension}
$$ \exists Z\, \forall x [x \in Z \leftrightarrow F(x)], $$
where $F(a)$ is an
{\it arithmetic formula}, i.e. a formula without set quantifiers.

 For a 2-place relation $\prec$ and an arbitrary formula $F(a)$ of $\mathcal{L}_2$ we define
\begin{enumerate}
\item[]
$\text{Prog}(\prec,F):=(\forall x)[\forall y (y\prec x \rightarrow F(y))\rightarrow F(x)]$ (\emph{progressiveness})
\item[]$\text{\bf{TI}}(\prec,F):= \text{Prog}(\prec,F)\rightarrow\forall x F(x)$ (\emph{transfinite induction})
\item[] $\text{WF}(\prec):=\forall X\text{\bf{TI}}(\prec,X):=$ \newline
$\forall X(\forall x[\forall y (y\prec x \rightarrow y\in X))\rightarrow x\in X]\rightarrow \forall x[x\in X])$ (\emph{well-foundedness}).
\end{enumerate}
Let $\mathcal{F}$ be any collection of formulae of $\mathcal{L}_2$. For a 2-place relation $\prec$ we will write $\prec\in\mathcal{F}$, if $\prec$   is defined by a formula $Q(x,y)$ of $\mathcal{F}$ via $x\prec y:=Q(x,y)$.

\begin{dfn}{\em
$\mathrm{BI}$ denotes the bar induction scheme, i.e. all formulae  of the form
$$\text{WF}(\prec)\rightarrow\text{\bf{TI}}(\prec,F),$$
where $\prec$ is an arithmetical relation (set parameters allowed) and $F$ is an arbitrary formula of $\mathcal{L}_2$.

By $\BI$ we shall refer to the theory $\ACA_0+\mathrm{BI}$.
}\end{dfn}

\begin{remark}\label{tarski-truth}{\em The statement of the main theorem \ref{main} uses the notion of a countable coded
$\omega$-model of $\BI$. As the stated equivalence is claimed to be provable in $\RCA_0$,
a few comments on how this is formalized in this weak base theory are in order.
The notion of a countable coded $\omega$-model can be formalized in $\RCA_0$ according to
 \cite[Definition VII.2.1]{SOSA}. Let $\mathbb{M}$ be a countable coded $\omega$-model.
 Since $\BI$ is not finitely axiomatizable we have to quantify over all axioms of $\BI$
 to express that $\mathbb{M}\models \BI$. The axioms of $\BI$ (or rather their G\"odel numbers)
 clearly form a primitive recursive set, $Ax(\BI)$.
 To express $\mathbb{M}\models \phi$ for $\phi\in Ax(\BI)$ we use the notion of a {\em valuation for} $\phi$ from \cite[Definition VII.2.1]{SOSA}.
 A valuation $f$ for $\phi$ is a function from the set of subformulae of $\phi$ into the set $\{0,1\}$ obeying the usual Tarski truth conditions.
 Thus we write $\mathbb{M}\models \phi$, if there exists a valuation $f$ for $\phi$ such
 that $f(\phi)=1$. Whence $\mathbb{M}\models \BI$ is defined by
 $\forall \phi\in Ax(\BI)\;\mathbb{M}\models \phi$.
 }\end{remark}

\section{Relativizing the Howard-Bachmann ordinal}\label{s2}
In this section we show how to relativize the construction that
leads to the Howard-Bachmann ordinal to an arbitrary countable
well-ordering. To begin with, mainly to foster intuitions, we provide a set-theoretic definition
working in $\mathbf{ZFC}$.
This will then be followed by a purely formal definition that can be made
in $\mathbf{RCA}_0$.

Throughout this section, we fix a countable well-ordering
$\mathfrak{X}=(X,<_X)$ without a maximum element, i.e.,   an ordered pair $\mathfrak{X}=(X,<_X)$, where $X$ is a set
of natural numbers, $<_X$ is  a well-ordering relation on $X$, and $\forall v\in X\,\exists u\in X\;v<_Xu$. We write $\dX$ for $X$.

Firstly, we need some ordinal-theoretic background. Let $\On$ be the class of
ordinals.
Let $\AP:=\{\xi\in \On\!:\exists\eta\in \On [\xi=\omega^{\eta}]\}$
be the class of additive principal numbers and
let $\Ep := \{ \xi \in \On\!:\xi = \omega^{\xi} \}$ be the class of
$\varepsilon$--numbers which is enumerated by the function $\lambda \xi .
\varepsilon_{\xi}$.

We write $\al =_{NF}\om^{\alpha_1}+\ldots+\om^{\alpha_n}$ if
$\al=\om^{\alpha_1}+\ldots+\om^{\alpha_n}$
and $\alpha>\alpha_1\geq \ldots\alpha_n$.
Note that by Cantor's normal form theorem, for every $\al \notin \Ep\cup \{ 0\}$,
there
are uniquely determined ordinals $\alpha_1,\ldots,\alpha_n$  such that
$\al =_{NF}\om^{\alpha_1}+\ldots+\om^{\alpha_n}$.

Let $\Omega := \aleph_{1}$. For $u\in \dX$, let $\EE_u$ be the $u^{th}$ $\varepsilon$-number $>\Omega$.
Thus, if $u_0$ is the smallest element of $\dX$, then $\EE_{u_0}$ is the least $\varepsilon$-number $>\Omega$, and in general,
for $u\in \dX$ with $u_0<_Xu$, $\EE_u$ is the least $\varepsilon$-number $\rho$ such that $\forall v<_Xu\;\EE_v<\rho$.

In what follows we shall only be interested in ordinals below $\sup_{u\in X}\EE_u$. Henceforth, unless indicated otherwise,
any ordinal will be assumed to be smaller than that ordinal.

For any such $\alpha$ we
define the set $E_{\Omega}(\alpha)$ which consists of the
$\varepsilon$--numbers below $\Omega$ which are needed for the unique
representation of $\alpha$ in Cantor normal form
 recursively as follows:
\begin{enumerate}
\item $E_{\Omega}(0) := E_{\Omega}(\Omega) := \emptyset $ and $E_{\Omega}(\EE_u):= \emptyset $ for $u\in \dX$.
\item $E_{\Omega}(\alpha) := \{ \alpha \}, \mbox{ if } \alpha \in E \cap
\Omega,$
\item $E_{\Omega}(\alpha) := E_{\Omega}(\alpha_1)
\cup\ldots \cup E_{\Omega}(\alpha_n)\;\;\mbox{ if }\;\;
\al =_{NF}\om^{\alpha_1}+\ldots+\om^{\alpha_n}$.
\end{enumerate}
Let $\alpha^{\astr} := \max(E_{\Omega}(\alpha) \cup \{ 0 \})$. \\[0.5cm]
We  define sets of ordinals $\CC(\alpha, \beta),
\CC^{n}(\alpha,\beta)$, and ordinals $\vt \alpha$ by main recursion
on $\alpha< \sup_{u\in X}\EE_u$ and subsidiary recursion on $n < \omega$
(for
$\beta
<
\Omega$) as follows.

\begin{enumerate}
\item [(C0)]  $\EE_u\in
\CC^{n}(\alpha,\beta)$ for all $u\in\dX$.
\item [(C1)]  $\{ 0,\Omega \} \cup \beta \;\subseteq\;
\CC^{n}(\alpha,\beta),$
\item [(C2)] $\gamma_1,\ldots,\gamma_n\in \CC^{n}(\alpha,\beta)\;\wedge\;\xi=_{NF}\om^{\ga_1}+\ldots+\om^{ga_n}
\;\Longrightarrow\;
\xi \in \CC^{n+1}(\alpha,\beta),$
\item [(C3)] $\delta \in \CC^{n}(\alpha,\beta) \cap \alpha \;\Longrightarrow\;
\vt \delta \in \CC^{n+1}(\alpha,\beta),$
\item [(C4)] $\CC(\alpha,\beta) := \bigcup \{ \CC^{n}(\alpha,\beta) \!: n <
\omega \} ,$
\item [(C5)] $\vt \alpha := \min \{ \xi < \Omega \!: \CC(\alpha,\xi)
\cap \Omega \subseteq \xi \;\wedge\; \alpha \in \CC(\alpha,\xi) \} $
if there exists an ordinal $\xi<\Omega$ such that
$\CC(\alpha,\xi)
\cap \Omega \subseteq \xi$ and  $\alpha \in \CC(\alpha,\xi)$.
Otherwise $\vt\al$ will be undefined.

We will shortly see that $\vt\al$ is always defined (Lemma \ref{1.2}).
\end{enumerate}

\begin{remark}{\em The definition of $\vt$ originated in \cite{ra89}. An ordinal representation
system based on $\vt$ was used in \cite{ra} to determine the proof-theoretic strength of fragments
of Kripke-Platek set theory and in \cite{rathjen-weiermann93}
it was used to characterize the strength of Kruskal's theorem.
}
\end{remark}
\begin{lemma}\label{1.2}
$\vt\alpha$ is defined for every $\alpha < \sup_{u\in X}\EE_u$.
\end{lemma}
{\sl Proof:} Let $\beta_{0} := \alpha^{\astr} + 1$. Then $\alpha \in \CC(\alpha,\beta_{0})$ via (C1) and
(C2). Since the cardinality of $\CC(\alpha,\beta)$ is less than
$\Omega$ there exists a $\;\beta_{1} < \Omega$ such that
$\CC(\alpha,\beta_{0}) \cap \Omega \subset \beta_{1}$. Similarly there
exists for each $\beta_{n} < \Omega$ (which is constructed recursively)
a $\beta_{n+1} < \Omega$ such that $\CC(\alpha,\beta_{n}) \cap \Omega
\subseteq \beta_{n+1}$. Let $\beta := \sup \{ \beta_{n} \!: n < \omega
\}$. Then $\alpha \in \CC(\alpha, \beta)$ and $\CC(\alpha,\beta) \cap \Omega
\subset \beta < \Omega$. Therefore $\vt  \alpha \le \beta <
\Omega$.\hfill $\Box$ \\
\begin{lemma}\label{1.3}
\begin{enumerate}
\item $\vt \alpha \in E,$
\item $\alpha \in \CC(\alpha, \vt \alpha),$
\item $\vt \alpha = \CC(\alpha, \vt \alpha) \cap \Omega,${\mbox
 { and }}
$\vt\al\notin \CC(\al,\vt\al)$,
\item $\ga\in \CC(\al,\be) \iff \ga^{\astr}\in \CC(\al,\be)$,
\item $\al^{\astr}<\vt\al$,
\item $\vt  \alpha = \vt  \beta \;\Longrightarrow\; \alpha =
\beta,$
\item $\begin{array}{lcl}\vt  \alpha < \vt  \beta & \iff&(\alpha <
\beta \;\wedge\; \alpha^{\astr} < \vt  \beta) \;\vee\; (\beta < \alpha
\;\wedge\; \vt  \alpha \le \beta^{\astr})
\\ &\iff &
(\alpha <
\beta \;\wedge\; \alpha^{\astr} < \vt  \beta) \;\vee\;  \vt  \alpha \le \beta^{\astr} \end{array}$
\item $\beta < \vt  \alpha \iff\omega^{\beta} <
\vt  \alpha.$
\end{enumerate}
\end{lemma}
{\sl Proof:} (1) and (8)  basically follow from closure of $\vt\al$ under
(C2).

(2)  follows from the definition of
$\vt\al$ taking Lemma \ref{1.2} into account.

For (3), notice that $\vt\al\subset \CC(\al,\vt\al)$ is a consequence of clause (C1). Since $\CC(\al,\vt\al)\cap\Omega\subseteq \vt\al$ follows from the definition of
$\vt\al$ and Lemma \ref{1.2}, we arrive at (3).

(4): If $\ga^{\astr}\in \CC(\al,\be)$, then $\ga\in \CC(\al,\be)$ by (C2).
On the other hand,\\  $\ga\in \CC^n(\al,\be)\;\Longrightarrow\;\ga^{\astr}\in \CC^n
(\al,\be)$ is easily seen by induction on $n$.

(5): $\al^{\astr}\in \CC(\al,\vt\al)$ holds by (4). As $\al^{\astr}
<\Om$, this implies $\al^{\astr}<\vt\al$ by (3).

(6): Suppose, aiming at a contradiction, that $\vt\al = \vt\be$ and
$\al <\be$. Then $\CC(\al,\vt\al) \;\subseteq\; \CC(\be,\vt\be)$; hence
$\al\in \CC(\be,\vt\be)\cap\be$ by (2); thence $\vt\al = \vt\be\in \CC
(\be,\vt\be)$, contradicting (3).

(7): Suppose $\al<\be$. Then $\vt\al <\vt\be$ implies $\al^{\astr}<\vt
\be$
by (5). If $\al^{\astr}<\vt\be$, then $\al\in \CC(\be,\vt\be)$; hence
$\vt\al\in \CC(\be,\vt\be)$; thus $\vt\al<\vt\be$. This shows
$$(a)\;\;\;\al<\be\;\Longrightarrow\;(\vt\al<\vt\be \iff \al^{\astr}<\vt\be).$$
By interchanging the roles of $\al$ and $\be$, and employing (6) (to exclude $\vt\al=\vt\be$), one
obtains
$$(b)\;\;\;\be<\al\;\Longrightarrow\;(\vt\al<\vt\be \iff \vt\al\le\be^{\astr}).$$
$(a)$ and $(b)$ yield the first equivalence of (7) and thus the direction
``$\Rightarrow$" of the second equivalence.
Since $\vt\al\le\be^{\astr}$ implies $\vt\al<\vt\be$ by (5), one also obtains the direction
``$\Leftarrow$" of the second equivalence.
\qed

\begin{definition}\label{1.4} {\em
Inductive definition of a set $\TT$ of ordinals and a natural number
$G_{\vt }\al$ for $\alpha\in \TT$.
\begin{enumerate}
\item $0,\Omega \in \TT,\;\Gt 0:=\Gt \Om:=0,$.
$\EE_u\in\TT$ and $\Gt\EE_u=0$
for all $u\in\dX$.
\item If $\alpha =_{NF}\omega^{\alpha_{1}} + \ldots + \omega^{\alpha_{n}}$ and $\alpha_1,\ldots,\alpha_n\in\TT$ then $\alpha\in\TT$
and $\Gt\alpha:=\max\{\Gt\alpha_1,\ldots,\Gt\alpha_n\}+1$.

\item If $\alpha = \vt  \alpha_{1}$ and $\alpha_{1} \in \TT$
then $\alpha \in
\TT$ and $\Gt\alpha:=\Gt\alpha_1+1.$
\end{enumerate}}
\end{definition}

Observe that according to Lemma \ref{1.3} (1) and \ref{1.3} (6) the function
$\Gt$ is well-defined. Each ordinal $\al \in \TT$ has a unique normal form
using the symbols $0,\Om,+,\om,\vt$.

\begin{lemma}\label{1.6}
$ \TT =\bigcup \{\CC(\al,0)\!:\al< \sup_{u\in X}\EE_u\}=\CC(\sup_{u\in X}\EE_u,0)$.
\end{lemma}
\prf Obviously $\beta<\sup_{u\in X}\EE_u$ holds for all $\beta\in\TT$.
$$\beta\in \TT \Rightarrow \beta\in\CC(\sup_{u\in X}\EE_u,0)$$ is then shown by induction on $\Gt\beta$.

The inclusion $\CC(\sup_{u\in X}\EE_u,0)\subseteq\TT$ follows from the fact
that $\TT$ is closed under the clauses (Ci) for $i=0,1,2,3$.
Since $\mathfrak X$ is an ordering without a maximal element
it is also clear that  $\bigcup \{\CC(\al,0)\!:\al< \sup_{u\in X}\EE_u\}=\CC(\sup_{u\in X}\EE_u,0)$. \qed

 If for $\al,\be\in\TT$
represented in their normal form,
we wanted to determine whether $\al<\be$, we could do this by deciding $\al_0<\be_0$ for
ordinals $\al_0$ and $\be_0$ that appear in these representations and, in
addition, satisfy $\Gt\al_0+\Gt\be_0<\Gt\al+\Gt\be$. This follows from
Lemma 1.2 (7) and the recursive procedure for comparing ordinals in
Cantor normal form. So we come to see that
after a straightforward coding in the natural numbers, we may
represent
$\langle \TT, < \restriction \TT \rangle$ via a primitive
recursive ordinal notation system.
How this ordinal representation system can be directly defined in $\mathbf{RCA}_0$
is spelled out  in the next subsection.

\subsection{Defining $\TT$ in $\mathbf{RCA}_0$}
We shall provide an explicit primitive recursive definition of $\TT$ as a term structure in $\mathbf{RCA}_0$. Of course formally, terms or strings of symbols have to be treated as coded by natural numbers since $\mathbf{RCA}_0$ only talks about
numbers and sets of numbers. Though, as it is well-known how to do this,
we can't be bothered with these niceties.

\begin{dfn}\label{2.2}
Given  a well-ordering
$\mathfrak{X}=(X,<_X)$, i.e.,   an ordered pair $\mathfrak{X}$ in  which $X$ is a set
of natural numbers
  and $<_X$ is  a well-ordering relation on $X$,
we define, by recursion, a binary relational structure  $\vtt=(\vert\vtt\vert,\kvt)$, and a function $\thickspace ^*:\vert\vtt\vert\to\vert\vtt \vert$,  in the following way:
\begin{enumerate}
\item  $ 0, \Omega \in \vert\vtt \vert$,
and $0^*:=0=:\Omega^*$.
\item If $\alpha\in\vd$ and $0\ne \alpha$ then $0\kvt \alpha$.
\item For every $u\in X$ there is an element $\EE_u\in\vert\vtt \vert$. Moreover, $(\EE_u)^*:=0$, and $\Omega\kvt \EE_u$.
If $u,v\in X$ and $u<_X v$, then $\EE_u<\EE_v$.

\item For every $\alpha\in\vert\vtt \vert$
there is an element $\vt\alpha\in\vert\vtt \vert$; and we have $\vt\alpha\kvt\Omega$, $\vt\alpha\kvt \EE_u$ for every $u\in X$, and $(\vt\alpha)^*:=\vt\alpha$.

\item If $\alpha\in \vd $ and $\alpha$ is not of the form $\Omega$, $\EE_u$, or $\vt\beta$, then $\omega^\alpha\in\vtt $ and
$(\omega^\alpha)^*:=\alpha^*$.

\item  If $\alpha_1,\ldots,\alpha_n\in \vd$ and $\alpha_1\kgrvt\dots\kgrvt\alpha_n$ with $n\geq 2$, then $\omega^{\alpha_1}+\omega^{\alpha_2}+\dots+\omega^{\alpha_n}\in\vd $ and
   $(\omega^{\alpha_1}+\omega^{\alpha_2}+\dots+\omega^{\alpha_n})^*:=\max\{\alpha_i^*:1\leq i\leq n\}$.
   \item Let $\alpha=\omega^{\alpha_1}+\dots+\omega^{\alpha_n}\in \vd$ and $\beta\in \vd$, where $\beta$ is of one of the forms
   $\vt \gamma$, $\Omega$, or $\EE_u$.
    \begin{itemize}
    \item[(i)] If $\alpha_1\kvt \beta$, then $\omega^{\alpha_1}+\dots+\omega^{\alpha_n}\kvt\beta$.
    \item[(ii)] If $\beta\kgvt \alpha_1$, then $\beta\kvt \omega^{\alpha_1}+\dots+\omega^{\alpha_n}$.
    \end{itemize}
\item If $\omega^{\alpha_1}+\dots+\omega^{\alpha_n}, \thickspace \omega^{\beta_1}+\dots+\omega^{\beta_m}\in\vert\vtt \vert$ then
\subitem $\omega^{\alpha_1}+\dots+\omega^{\alpha_n}\kvt\omega^{\beta_1}+\dots+\omega^{\beta_m}$ iff
\subitem $n<m\medspace\land\medspace\forall i\leq n \thickspace\thinspace \alpha_i=\beta_i$ or
\subitem $\exists\thinspace i\leq\min(n,m) \thickspace [(\forall j<i \thickspace \alpha_j=\beta_j)\land(\alpha_i<\beta_i)]$.
\item If $\alpha\kvt\beta$ and $\alpha^*\kvt\vt\beta$ then $\vt\alpha\kvt\vt\beta$.
\item If $\vt\beta\kgvt \alpha^*$ then $\vt\beta\kvt\vt\alpha$.

\end{enumerate}
\end{dfn}

\begin{lemma} \begin{itemize}
\item[(i)] The set $\vd$, the relation $\kvt$, and the function $^*$ are primitive recursive in $\mathfrak{X}=(X,<_X)$.
\item[(ii)] $\kvt$ is a total and linear ordering on $\vd$.
\end{itemize}
\end{lemma}
\proof Straightforward but tedious. \qed
Of course, $\mathbf{RCA}_0$ does not prove that
$\kvt$ is a well-ordering on $\vd$.

\section{A Well-ordering Proof}\label{s3}
In this section we work in the background theory  $$\RCA_0+\forall X\exists Y \;(X \in Y \und \mbox{$Y$ is an $\omega$-model of $\BI$})$$ and  shall prove the following statement
$$\forall {\mathfrak X}\,(\WO({\mathfrak X})\to \WO(\vt_{\mathfrak X}))\,,$$
that is, the part $\mbox{(i)}\Rightarrow \mbox{(ii)}$ of the main theorem \ref{main}. Some of the proofs are similar to ones in \cite{rathjen-weiermann93} section 10.
Note that in this theory we can deduce arithmetical comprehension and even arithmetical transfinite recursion
owing to  \cite{girard87} and \cite{H-Band}, respectively.

Let us fix a well-ordering ${\mathfrak X}=(X,<_X)$, an arbitrary set $Y$ and a countable coded $\omega$-model $\Mod$ of $\BI$ which contains both ${\mathfrak X}$ and $Y$
as  elements.
In the sequel $\al,\be,\ga,\de,\ldots$ are supposed to range over $\vt_{\mathfrak X}$.
$<$ will be used to denote the ordering on $\vt_{\mathfrak X}$.
We are going to work informally in our background theory.
A set $U\subseteq {\mathbb N}$ is said to be definable in $\Mod$ if $U=\{n\in{\mathbb N}\mid {\Mod}\models A(n)\}$
for some formula $A(x)$ of second order arithmetic which may contain  parameters from $\Mod$.

\begin{definition}{\em
\begin{enumerate}
\item$\Acc:=\{\al<\Om\mid {\Mod}\models \WO(<\restriction \al)\},$
\item$\Mm:=\{\al\!:E_{\Om}(\al)\subseteq \Acc\},$
\item$\al<_{\Om}\be:\iff\al,\be\in \Mm\;\wedge\;\al<\be.$
\end{enumerate}}
\end{definition}

\begin{lemma}\label{5.2}
$\al,\be\in \Acc\;\Longrightarrow\;\al+\om^{\be}\in \Acc.$
\end{lemma}
{\em Proof.} Familiar from Gentzen's proof in Peano arithmetic. The proof just requires $\ACA_0$.
(cf. \cite[VIII.\textsection 21 Lemma 1]{sc}).\qed
\begin{lemma}\label{10.2}
$\Acc =\Mm\cap\Om\; (:=\{\al\in \Mm\mid \al<\Om\}.)$
\end{lemma}
{\em Proof.}
If $\al\in \Acc$, then $E_{\Om}(\al)\subseteq \Acc$
as well; hence $\al\in \Mm\cap \Om$. If $\al\in \Mm\cap \Om$, then $E_{\Om}(\al)\subseteq
\Mm\cap \Om$, so $\al\in Acc $ follows from Lemma \ref{5.2}. \qed

\begin{lemma}\label{mod-bi} Let $U$ be $\Mod$ definable. Then
$$\forall \al<\Omega\cap \Mm \,[\forall \beta<\alpha \beta\in U \to \alpha\in U]\to \Acc\subseteq U\,.$$
\end{lemma}
\proof This follows readily from the assumption that $\Mod$ is a model of $\BI$. \qed

\begin{definition}{\em
Let $\Prog_{\Om}(X)$ stand for
$$(\forall \al\in \Mm)[(\forall \be\kOm\al) (\be\in X)\longrightarrow\al\in X].$$
Let $\Acc_{\Om}:=\{\al\in \Mm\!:\vt\al\in \Acc\}.$}
\end{definition}

\begin{lemma}\label{Omega} If $U$ is $\Mod$ definable, then
$$\Prog_{\Omega}(U)\to \Omega,\Omega+1\in U\,.$$
\end{lemma}
\prf This follows from Lemma \ref{10.2} and Lemma \ref{mod-bi}. \qed

\begin{lemma}\label{10.3}
$\Prog_{\Om}(\Acc_{\Om}).$
\end{lemma}
\prf
Assume $\al\in \Mm$ and $(\forall \be\kOm\al)(\be\in \Acc_{\Om}).$
We have to show that $\vt\al\in \Acc.$ It suffices to show
\begin{eqnarray}
\label{we1}
\be<\vt\al &\Longrightarrow &\be \in \Acc.
\end{eqnarray}
We shall employ induction on $\Gt(\be)$, i.e.,
 the length of (the term that represents) $\be$.
If $\be\not\in E$, then (\ref{we1}) follows easily by the inductive assumption and
Lemma
\ref{5.2}. Now suppose $\be=\vt\be_0.$
According to Lemma \ref{1.3} it suffices to consider the following two cases:\\
{\em Case 1:} $\be\le \al^*$. Since $\al\in \Mm,$ we have $\al^*\in E_{\Om}(\al)
\subseteq \Acc;$ therefore $\be\in \Acc.$\\
{\em Case2:} $\be_0<\al$ and $\be_0^*<\vt\al.$
As the length of $\be_0^*$ is less than the length of $\be$, we get $\be_0^*\in \Acc;$
thus $E_{\Om}(\be_0)\subseteq \Acc,$ therefore $\be_0\in \Mm.$ By the assumption
at the beginning of the proof, we then get $\be_0\in \Acc_{\Om};$ hence
$\be=\vt\be_0\in \Acc$.\qed

\begin{definition}\label{Fj}{\em For every $\Mod$ definable set $U$ we define the ``Gentzen  jump''
 $$U^j \;:=\; \{\gamma\mid \forall \delta\,[\Mm\cap\delta\subseteq  U\to
\Mm\cap(\delta+\omega^\gamma)\subseteq U ]\}.$$
}\end{definition}

\begin{lemma}\label{12.6} Let $U$ be $\Mod$ definable.
\begin{itemize}
\item[(i)] $\gamma\in U^j\Rightarrow \Mm\cap\omega^\gamma\subseteq  U$.
\item[(ii)] $\Prog_{\Omega}(U)\Rightarrow \Prog_{\Omega}(U^j)$.
\end{itemize}
\end{lemma}
 \prf (i) is obvious. (ii)  $\Mm\cap(\delta+\omega^\gamma)\subseteq U$ is to be proved under the assumptions (a) $\Prog_{\Omega}(U)$, (b) $\gamma\in \Mm\und \Mm\cap\gamma\subseteq U^j$ and (c) $\Mm\cap \delta\subseteq U$. So let $\eta\in \Mm\cap (\delta+\omega^{\gamma})$.
\begin{enumerate}
\item $\eta<\delta$: Then $\eta\in U$ is a consequence of (c).
\item $\eta=\delta$: Then $\eta\in U$ follows from (c) and (a).
\item $\delta<\eta<\delta+\omega^\gamma$: Then there exist $\gamma_1,\ldots,\gamma_k<\gamma$ such that $\eta=\delta+\omega^{\gamma_1}+\ldots+\omega^{\gamma_k}$ and $\gamma_1\geq \ldots\geq\gamma_k$. $\eta\in \Mm$ implies $\gamma_1,\ldots,\gamma_k\in \Mm\cap\gamma$. Through applying (b) and (c) we obtain $\Mm\cap(\delta+\omega^{\gamma_1})\subseteq U$. By iterating this procedure we  eventually arrive at $\delta+\omega^{\gamma_1}+\ldots+\omega^{\gamma_k}\in U$, so $\eta\in U$ holds. \end{enumerate}\qed

\begin{cor}\label{12.7} Let $\Ii(\delta)$ be the statement that  $\Prog_{\Omega}(V)\to \delta\in \Mm\und \delta\cap \Mm\subseteq V$ holds for all $\Mod$ definable sets $V$.
Assume $\Ii(\delta)$. Let $\delta_0:=\delta$ and $\delta_{n+1}:=\omega^{\delta_n}$.
Then $$\Ii(\delta_n)$$
holds for all $n$.
\end{cor}
\prf We use induction on $n$. For $n=0$ this is the assumption. Now suppose $\Ii(\delta_n)$ holds.
Assume $\Prog_{\Omega}(U)$ for an $\Mod$ definable $U$. By Lemma \ref{12.6} we conclude $\Prog_{\Omega}(U^j)$
and hence $\delta_n\in U^j$ and $\delta_n\cap \Mm \subseteq U^j$.
As clearly $\Mm\cap 0\subseteq U$ we get $\omega^{\delta_n}\cap \Mm \subseteq U$. Since $\Prog_{\Omega}(U)$ entails  $\delta\in\Mm$ we also have $\delta_{n+1}\in \Mm$. Thus $\delta_{n+1}\in \Mm\und \delta_{n+1}\cap \Mm \subseteq U$, showing
$\Ii(\delta_{n+1})$. \qed

Let $\omega_0(\alpha):=\alpha$ and $\omega_{n+1}(\alpha):=\omega^{\omega_n(\alpha)}$.

\begin{proposition}\label{Ii}  $\Ii(\EE_u)$ holds for all $u\in \vert {\mathfrak X}\vert$.
\end{proposition}
\prf Noting that in our background theory $\mathfrak X$ is a well-ordering, we can use induction on $\mathfrak X$.
Note also that $\Ii(\EE_u)$ is a statement about all definable sets in $\Mod$ which is not formalizable in $\Mod$ itself.
However, in our background theory quantification over all these sets is first order expressible and therefore transfinite
induction along $<_X$ is available.

First observe that we have $\Ii(\Omega+1)$ by Lemma \ref{Omega}. Let $u_0$ be the $<_X$-least element of $\vert {\mathfrak X}\vert$.
We have $\EE_{u_0}\in \Mm$ and for every $\eta<\EE_{u_0}$ there exists $n$ such that $\eta<\omega_n(\Omega+1)$.
As a result, using Corollary \ref{12.7}, we have $$\Prog_{\Omega}(U)\to \EE_{u_0}\cap \Mm\subseteq U$$ for every
$\Mod$ definable set $U$.

Now suppose that $u\in \vert {\mathfrak X}\vert$ is not the $<_X$-least element and for all $v<_X u$ we have
 $\Ii(\EE_v)$. As for every $\delta<\EE_u$ there exists $v<_X u$ and $n$ such that $\delta <\omega_n(\EE_v)$, the inductive
 assumption together with Corollary \ref{12.7} yields
 $$\Prog_{\Omega}(U)\to \EE_{u}\cap \Mm\subseteq U\,.$$ $\EE_u\in\Mm$ is obvious.
 \qed

 \begin{prop} For all $\alpha$, $\Ii(\alpha)$.
 \end{prop}
 \prf We proceed by the induction on the term complexity of $\alpha$.
 Clearly, $\Ii(0)$. By Lemma \ref{Omega} we conclude that   $\Ii(\Omega)$.
 Proposition \ref{Ii} entails
  that $\Ii(\EE_u)$ for all $u\in \vert{\mathfrak X}\vert$.

 Now let $\alpha=\omega^{\alpha_1}+\cdots+\omega^{\alpha_n}$ be in Cantor normal form. Inductively we have
 $\Ii(\alpha_1),\ldots,\Ii(\alpha_n)$. Assume $\Prog_{\Omega}(U)$. Then $\Prog_{\Omega}(U^j)$ by Lemma \ref{12.6}(ii),and hence
 $\alpha_1\cap M\subseteq U^j,\ldots,
 \alpha_n\cap M\subseteq U^j$ and $\alpha_1,\ldots, \alpha_n\in \Mm$. The latter implies $\alpha_1\in U^j,\ldots,
 \alpha_n\in U^j$. Using the definition of $U^j$ repeatedly we conclude
 $\alpha\cap M\subseteq U$. Moreover, $\alpha\in\Mm$ since  $\alpha_1,\ldots, \alpha_n\in \Mm$.

Now suppose that $\alpha=\vt \beta$. Inductively we have $\Ii(\beta)$. By Lemma \ref{10.3} we conclude that $\beta\in \Acc_{\Omega}$,
and hence $\alpha\in \Acc$. From $\Prog_{\Omega}(U)$ we obtain by Lemma \ref{mod-bi} that $\xi\in U$ for all
$\xi\leq\alpha$. As a result, $\Ii(\alpha)$. \qed

\begin{cor} $\vt_{\mathfrak X}$ is a well-ordering.
\end{cor}
With the previous Corollary, the proof of Theorem \ref{main} (i)$\Rightarrow$(ii) is finally accomplished.

\section{Deduction chains}\label{s4}
From now on we will be concerned with the part $\mbox{(ii)}\Rightarrow \mbox{(i)}$ of the main theorem \ref{main}.
An important tool will be the method of deduction chains.
Given a sequent $\Gamma$  and a set $Q\subseteq \mathbb N$, deduction chains
starting at $\Gamma$  are built by systematically decomposing  $\Gamma$ into its subformulas, and adding additionally at the
$n$th step the formulas $\neg A_n$ and $\neg \bar{Q}(\bar{n})$, where $(A_n \mid n\in \mathbb N)$ is an enumeration of the
axioms of the theory $\BI$, and $\bar{Q}(\bar{n})$ is the atom $\bar{n}\in U_0$ if $n \in Q$ and $\bar{n}\notin U_0$
otherwise. The set of all deduction chains that can be built from the empty sequent
with respect to  a given set $Q$ forms the tree ${\mathcal D}_Q$. There are two
scenarios to be considered.
\begin{itemize}
\item[(i)] If there is an infinite deduction chain, i.e. ${\mathcal D}_Q$ is ill-founded, then this readily
yields a model of $\BI$ that contains $Q$.
\item[(ii)] If each deduction chain is finite, then this yields a derivation of the empty sequent, $\perp$, in a corresponding
infinitary system with an $\omega$-rule. The depth of this derivation is
bounded by the order-type $\alpha$ of the Kleene-Brouwer ordering of ${\mathcal D}_Q$. By the
well-ordering principle, transfinite induction up to ${\mathfrak E}_{\alpha+1}$ is available, which allows
to transform this proof into a cut-free proof of $\perp$ whose depth is less
than $\vartheta{\mathfrak E}_{\alpha+1}$.
\end{itemize}
As the second alternative is impossible, the first yields the desired model.
\begin{dfn}{\em
\begin{enumerate}
\item We let $U_0, U_1, \ldots,
U_m,\ldots
$ be an enumeration of the free set variables of $\mathcal{L}_2$ and, given a closed term $t$, we write $t^{\nel}$ for its numerical value.
\item Henceforth a {\bf{sequent}} will be a finite set of $\mathcal{L}_2$-formulae \emph{without} free number variables.
\item A sequent $\Gamma$ is {\bf{axiomatic}} if it satisfies at least one of the following conditions:
\begin{enumerate}
\item $\Gamma$ contains a true {\bf{literal}}, i.e., a true formula of either of the  forms $R(t_1,\dots,t_n)$ or $\neg R(t_1,\dots,t_n)$, where $R$ is a predicate symbol in $\mathcal{L}_2$ for a primitive recursive relation and $t_1,\dots, t_n$ are closed terms.
\item $\Gamma$ contains the  formulae $s\in U$ and
$t \notin U$ for some set variable $U$  and terms $s,\medspace t$ with $s^{\nel}=t^{\nel}$.
\end{enumerate}
\item A sequent is {\bf{reducible}}
if it is not axiomatic and contains a formula which is not a literal.
\end{enumerate}
}\end{dfn}
\begin{dfn}{\em
For $Q\subseteq \nel$ we define
\begin{equation*}
\bar Q(n) \leftrightharpoons
\begin{cases}
\bar n\in U_0 & \text{if } n\in Q,\\
\bar n\notin U_0 & \text{otherwise}
\end{cases}
\end{equation*}
}\end{dfn}

For some of the following theorems it is convenient to have a finite axiomatization of arithmetical comprehension.

\begin{lemma}
$\ACA_0$ can be axiomatized via a single $\Pi^1_2$ sentence $\forall X C(X)$.
\end{lemma}
\prf \cite[Lemma VIII.1.5]{SOSA}. \qed

\begin{dfn}
In what follows, we fix an enumeration of $A_1, \thinspace A_2,\thinspace A_3, \dots$ of all the universal closures of instances of {\rm({\bf{BI}})}. We also put $A_0:=\forall X \thinspace C(X)$, where the latter is the sentence  axiomatizes arithmetical comprehension.
\end{dfn}

\begin{dfn}{\em
Let $Q\subseteq \nel$.
A  $Q${\bf{-deduction chain}} is a finite string
$$\Gamma_0,\thinspace \Gamma_1,\dots,\thinspace \Gamma_k$$
of sequents $\Gamma_i$ constructed according to the following rules:
\begin{enumerate}
\item $\Gamma_0=\neg\bar Q(0),\thinspace \neg A_0$. 
\item $\Gamma_i$ is not axiomatic for $i<k$.
\item If $i<k$ and $\Gamma_i$ is not reducible then
$$\Gamma_{i+1}=\Gamma_i,\thinspace \neg\bar Q(i+1),\thinspace \neg A_{i+1}$$ 
\item Every reducible $\Gamma_i$ with $i<k$ is of the form
$$\Gamma_i', \thinspace E,\thinspace \Gamma_i''$$
where $E$ is not a literal and $\Gamma_i'$ contains only literals.
$E$ is said to be the {\bf{redex}} of $\Gamma_i$.

Let $i<k$ and $\Gamma_i$ be reducible.
$\Gamma_{i+1}$ is obtained from $\Gamma_i=\Gamma_i',\thinspace E,\thinspace \Gamma_i''$ as follows:
\begin{enumerate}
\item If $E\equiv E_0\lor E_1$ then
$$\Gamma_{i+1}=\Gamma_i',\thinspace E_0,\thinspace E_1,\thinspace \Gamma_i'',\thinspace \neg \bar Q(i+1),\thinspace\neg A_{i+1}.$$
\item If $E\equiv E_0\land E_1$ then
$$\Gamma_{i+1}=\Gamma_i',\thinspace E_j,\thinspace \Gamma_i'',\thinspace \neg \bar Q(i+1),\thinspace\neg A_{i+1}$$
where $j=0$ or $j=1$.
\item If $E\equiv \exists x F(x)$ then
$$\Gamma_{i+1}=\Gamma_i',\thinspace F(\bar m),\thinspace \Gamma_i'',\thinspace \neg \bar Q(i+1),\thinspace \neg A_{i+1},\thinspace E$$
where $m$ is the first number such that $F(\bar m)$ does not occur in $\Gamma_0,\dots,\thinspace\Gamma_i$.
\item If $E\equiv \forall x F(x)$ then
$$\Gamma_{i+1}=\Gamma_i',\thinspace F(\bar m),\thinspace \Gamma_i'',\thinspace \neg \bar Q(i+1),\thinspace \neg A_{i+1}$$
for some $m$.
\item If $E\equiv\exists X F(X)$ then
$$\Gamma_{i+1}=\Gamma_i',\thinspace F(U_m),\thinspace
\Gamma_i'', \thinspace \neg \bar Q(i+1),\thinspace \neg A_{i+1},\thinspace E$$
where $m$ is the first number such that $F(U_m)$ does not occur in $\Gamma_0,\dots,\thinspace \Gamma_i$.

\item  If  $E \equiv \forall  X F(X)$  then
$$\Gamma_{i+1}=\Gamma_i',\thinspace F(U_m),\thinspace \Gamma_i'',\thinspace \neg \bar Q(i+1),\thinspace \neg A_{i+1}$$
where $m$ is the first number such that 
$U_m$ does not occur in $\Gamma_i$.
\end{enumerate}
\end{enumerate}
}\end{dfn}
The set of $Q$-deduction chains forms a tree $\mathcal{D}_Q$ labeled with strings of sequents.

We will now consider two cases.
\\[2ex] {\bf Case I:} ${\mathcal D}_Q$ is not well-founded.
 Then ${\mathcal D}_Q$ contains an infinite path $\mathbb P$. Now define
 a set $M$ via
 \begin{eqnarray*} (M)_i &=& \{k\mid\mbox{$\bar{k}\notin U_i$ occurs in ${\mathbb P}$}\}.\end{eqnarray*}
 Set ${\mathbb M}=({\mathbb N};\{(M)_i\mid i\in {\mathbb N}\},\in,+,\cdot,0,1,<)$.

For a formula $F$,
 let $F\in {\mathbb P}$ mean that $F$ occurs in $\mathbb P$, i.e. $F\in \Gamma$ for some $\Gamma \in {\mathbb P}$.
\\[1ex] {\bf Claim:} Under the assignment $U_i\mapsto (M)_i$ we have
 \begin{eqnarray}\label{claim} F\in {\mathbb P} &\;\;\;\;\Rightarrow\;\;\;\;& {\mathbb M}\models
 \neg F.\end{eqnarray}
 The Claim will imply that ${\mathbb M}$ is an $\omega$-model of
$\mathbf{BI}$. Also note that $(M)_0=Q$, thus $Q$ is in ${\mathbb M}$.
The proof of (\ref{claim}) follows by induction on $F$ using Lemma
 \ref{nextb} below.
The upshot of the foregoing is that we can prove Theorem \ref{main}
under the assumption that ${\mathcal D}_Q$ is ill-founded for all sets
$Q\subseteq{\mathbb N}$.

\begin{lemma}\label{nextb}
Let $Q$ be an arbitrary subset of $\bbM{N}$ and $\mathcal{D}_Q$ be the corresponding deduction tree. Moreover, suppose $\mathcal{D}_Q$ is not well-founded. Then $\mathcal{D}_Q$ has an infinite path $\bbM{P}$. $\bbM{P}$ has the following properties:
\begin{enumerate}
\item $\bbM{P}$ does not contain literals which are true in $\bbM{N}$.
\item $\bbM{P}$ does not contain formulas $s\in U_i$ and $t\notin U_i$ for constant terms $s$ and $t$ such that $s^\bbM{N}=t^\bbM{N}$.
\item If $\bbM{P}$ contains $E_0\lor E_1$ then $\bbM{P}$ contains $E_0$ and $E_1$.
\item If $\bbM{P}$ contains $E_0\land E_1$ then $\bbM{P}$ contains $E_0$ or $E_1$.
\item If $\bbM{P}$ contains $\exists xF(x)$ then $\bbM{P}$ contains $F(\bar{n})$ for all $n$.
\item If $\bbM{P}$ contains $\forall xF(x)$ then $\bbM{P}$ contains $F(\bar{n})$ for some  $n$.
\item If $\bbM{P}$ contains $\exists XF(X)$ then $\bbM{P}$ contains $F(U_m)$ for all  $m$.
\item If $\bbM{P}$ contains $\forall XF(X)$ then $\bbM{P}$ contains $F(U_m)$ for some  $m$.
\item $\bbM{P}$ contains $\neg C(U_m)$ for all  $m$.
\item $\bbM{P}$ contains $\neg \DQ(m)$ for all  $m$.
\end{enumerate}
\end{lemma}
\prf Standard. \qed

\begin{cor}\label{coro} If $\mathcal{D}_Q$ is ill-founded then there exists a
countable coded $\omega$-model of $\mathbf{BI}$ which contains $Q$.
\end{cor}

For our purposes it is important that Corollary \ref{coro} can be proved in $T_0:=\mathbf{RCA}_0+\forall {\mathfrak X}\,(\WO({\mathfrak X})\to \WO(\vt_{\mathfrak X}))$. To this end we need to show that the semantics
of $\omega$-models can be handled in the latter theory, i.e. for every
formula $F$ of $\mathcal{L}_2$ there exists a valuation for $F$
in the sense of \cite[VII.2.1]{SOSA}.
It is easily seen that the principle $\forall {\mathfrak X}\,(\WO({\mathfrak X})\to \WO(\vt_{\mathfrak X}))$ implies
$$\forall {\mathfrak X}\,(\WO({\mathfrak X})\to \WO(\varepsilon_{\mathfrak X}))$$ (see \cite[Definition 2.1]{AR}) and thus, by  \cite[Theorem 4.1]{AR}, $T_0$ proves that
every set is contained in an $\omega$-model of $\mathbf{ACA}$.
Now take an $\omega$-model containing $\mathcal{D}_Q$ and an infinite branch of  $\mathcal{D}_Q$. In this $\omega$-model
we find a valuation for every formula by \cite[VII.2.2]{SOSA}.
And hence Corollary \ref{coro} holds in the model, but then it also holds
in the world at large by absoluteness.

\section{Proof of the Main Theorem: The hard direction part 2}\label{s5}
The remainder of the paper will be devoted to ruling out the possibility that for some $Q$, ${\mathcal D}_Q$ could be a well-founded tree.
This is the place where the principle $\forall {\mathfrak X}\,(\WO({\mathfrak X})\to \WO(\vt_{\mathfrak X}))$
 in the guise of cut elimination for an
infinitary proof system  enters the stage.
Aiming at a contradiction, suppose that ${\mathcal D}_Q$ is a well-founded tree.
Let $\fXO$ be the Kleene-Brouwer ordering on  ${\mathcal D}_Q$ (see \cite[Definition V.1.2]{SOSA}).
Then $\fXO$ is a well-ordering.
In a nutshell, the idea is that a  well-founded ${\mathcal D}_Q$ gives rise
 to a derivation of the empty sequent (contradiction) in an infinitary proof system.

\subsection{Majorization and Fundamental Functions}
In this section we introduce the concepts of majorization
and fundamental function.
They are needed for carrying through the ordinal analysis of
bar induction. More details can be found in \cite{rathjen-weiermann93} section 4
and \cite[I.4]{bu} to which we refer for proofs.
The missing proofs are actually straightforward consequences of Definition
\ref{2.2}.
\begin{definition}{\em
\begin{enumerate}
\item$\alpha\lhd\beta$ means $\alpha<\beta$ and $\vt \alpha<\vt\beta$.
\item$\alpha\unlhd\beta:\;\Longleftrightarrow\;(\alpha\lhd\beta\vee
\alpha=\beta).$
\end{enumerate}}
\end{definition}
\begin{lemma}
\begin{enumerate}
\item$\alpha\lhd \beta\und \beta\lhd \gamma\;\Longrightarrow\;\alpha\lhd\gamma$.
\item$0<\beta<\varepsilon_0\;\Longrightarrow\;\alpha\lhd\alpha+\beta$.
\item$\alpha<\beta<\Omega\;\Longrightarrow\;\alpha\lhd\beta$.
\item$\alpha\lhd\beta\;\Longrightarrow\;\alpha+1\unlhd\beta.$
\item $\alpha\lhd\beta\;\Longrightarrow\;\vt\alpha\lhd\vt\beta.$
\item $\alpha=\alpha_0+1\Longrightarrow\vt\alpha_0\lhd\vt\alpha.$
\end{enumerate}
\end{lemma}
\
\begin{lemma}\label{4.4}
$\alpha\lhd\beta,\;\beta<\omega^{\gamma+1}\;\Longrightarrow\;\omega^{\gamma}+\alpha\lhd\omega^{\gamma}+\beta.$
\end{lemma}
\begin{corollary}
$\omega^{\alpha}\cdot n\lhd\omega^{\alpha}\cdot(n+1).$
\end{corollary}
\begin{lemma}\label{4.5}
$\alpha\lhd\beta\;\Longrightarrow\;\omega^{\alpha}\cdot
n\lhd\omega^{\beta}.$
\end{lemma}

\begin{definition}{\em Let $\AO:= (\TT\cap\Omega)\,\cup\,\{\Omega\}$. A function $f:\AO\to \TT$ will be called a {\em fundamental function} if it is generated by the following clauses:

\begin{enumerate}
\item[F1.] $Id:\AO\to\AO$ with $Id(\al)=\al$ is a fundamental function.
\item[F2.] If $f$ is a fundamental function, $\gamma\in\TT$ and $f(\Omega)<\omega^{\gamma+1}$,  then  $\omega^{\gamma}+f$ is a fundamental function, where
$(\omega^{\gamma}+f)(\alpha):=\omega^{\gamma}+f(\alpha)$ for all
$\alpha\in \AO.$
\item[F3.] If $f$ is a fundamental function then so is $\omega^{f}$ with $(\omega^{f})(\alpha):=\omega^{f(\alpha)}$
for all $\alpha\in \AO$.
\end{enumerate}}
\end{definition}
\begin{lemma}\label{6.7}
Let $f$ be a fundamental function and $\be\leq \Omega$.
\begin{itemize} \item[(i)] If $\alpha<\be$, then $f(\al)< f(\beta)$.
\item[(ii)] If $\alpha\lhd\be$, then $f(\al)\lhd f(\beta)$.
\item[(iii)] $(f(\beta))^{\astr}\leq \max((f(0))^{\astr},\beta^{\astr})$.
\end{itemize}
\end{lemma}
\proof (i) is obvious by induction on the generation of fundamental functions.

(ii) also follows by induction on the generation of fundamental functions, using Lemmata \ref{4.4} and \ref{4.5}.

(iii) as well follows by induction on the generation of fundamental functions.
\qed

\begin{lemma}\label{f5.4}
For every fundamental function $f$ we have
$f(\vt(f(0)))\lhd f(\Om)$.
\end{lemma}

\proof Since $\vt(f(0))<\Om$, we clearly have $f(\vt(f(0)))<f(\Om)$. Since $0\lhd \Om$ and $f$ is a fundamental function, we have $\vt(f(0))< \vt(f(\Omega))$ by lemma \ref{6.7} (ii). Invoking Lemma \ref{6.7} (iii), the latter entails that $(f(\vt(f(0))))^*< \vt(f(\Omega))$, so that in conjunction with  $f(\vt(f(0)))< f(\Om)$ it follows
that $\vt(f(\vt(f(0))))\lhd \vt(f(\Om))$. As a result, $f(\vt(f(0)))\lhd f(\Om)$.\qed

\subsection{The infinitary calculus $T^*_{\!_Q}$}
The calculus $T^*_{\!_Q}$ to be introduced stems from \cite{rathjen-weiermann93}
section 6.
We fix a set $Q\subseteq \mathbb N$.
Let ${\cal L}_2^Q$ be the language of second order arithmetic augmented by a unary predicate $\bar{Q}$.
The {\em formulas} of $T^*_{\!_Q}$ arise from ${\cal L}_2^Q$-formulas by
replacing
free numerical variables by numerals, i. e. terms of the form $0,0',0'',...$
Especially, every formula $A$ of $T^*_{\!_Q}$ is an
${\cal L}_2^Q$-formula.  We are going to measure the
length of derivations by ordinals.
We are going to use the set of ordinals $\TT$ of Section 3.
\begin{definition}{\em
\begin{enumerate}
\item A formula $B$ is said to be {\em weak} if it belongs to
$\Pi^1_0\cup\Pi^1_1$.
\item Two closed terms $s$ and $t$ are said to be equivalent if they
yield the same value when computed.
\item A formula is called constant if it contains no set variables.
The truth or falsity of such a formula is understood with respect to the
standard structure of the integers.
\item $\overline{0}:=0$, $\overline{m+1}:=\overline{m}'$.
\end{enumerate}}
\end{definition}
In the sequent calculus $T^*_{\!_Q}$ below we shall use the following rules of inference:\\[2ex]
$\begin{array}{ll}
(\wedge)&\vdash \Gamma,A$ and $\vdash \Gamma ,B
\;\Longrightarrow\; \vdash \Gamma, A \wedge B,\\[0.2cm]
(\vee)&\vdash \Gamma,A_{i} \;\Longrightarrow\; \vdash\Gamma, A_{0} \vee
A_{1} \;\;\; \mbox{ if }i \in \{ 0,1 \},\\[0.2cm]
(\forall_{2})&\vdash \Gamma, F(U) \;\Longrightarrow\; \vdash\Gamma,
\forall X F(X),\\[0.2cm]
(\exists_{1})&\vdash \Gamma, F(t) \;\Longrightarrow\; \vdash \Gamma,
\exists x F(x),\\[0.2cm]
(Cut)& \vdash \Gamma,A \mbox{ and } \vdash \Gamma, \neg \;A
\;\Longrightarrow\; \vdash \Gamma,
\end{array}$ \\[0.2cm]
where in  $(\forall_{2})$ the free variable  $U$
is not to occur in the conclusion.\\[1ex]
The most important feature of sequent calculi is cut--elimination.  To state this fact concisely,
let us introduce a measure of complexity, $gr(A)$, the {\em grade of a
formula A}, for ${\cal L}_2^Q$-formulae.

\begin{definition}\label{8.2} {\em
\begin{enumerate}
\item [1.] $gr(A) = 0$ if $A$ is a prime formula or negated prime
formula.
\item [2.] $gr(\forall X F(X)) = gr(\exists X F(X)) = \omega$ if $F(U)$
is arithmetic.
\item [3.] $gr(A \wedge B) = gr(A \vee B) = max \{ gr(A), gr(B) \} +1$.
\item [4.] $gr(\forall x H(x)) = gr(\exists x H(x)) = gr(H(0)) + 1$.
\item [5.] $gr(\forall X G(X)) = gr(\exists X G(X)) = gr(G(U)) + 1,$
\\ if $G$ is not arithmetic.
\end{enumerate}
}\end{definition}

\begin{definition}\label{8.3}{\em
Inductive definition of $\Tsr{\alpha}{\Gamma}$ for $\alpha\in\TT$ and
$\varrho<\omega+\omega$.\\
\begin{enumerate}
\item If $A$ is a true constant prime formula or negated prime formula and
$A\in \Gamma$, then $\Tsr{\alpha}{\Gamma.}$\\
\item If $n\in Q$ and $t$ is a closed term with value $n$ and $\bar{Q}(t)$ is in $\Gamma$, then
 $\Tsr{\alpha}{\Gamma.}$\\
 \item If $n\notin Q$ and $t$ is a closed term with value $n$ and $\neg \bar{Q}(t)$ is in $\Gamma$, then
 $\Tsr{\alpha}{\Gamma.}$\\
\item If $\Gamma$ contains formulas $A(s_1,\ldots,s_n)$ and $\neg
A(t_1,\ldots,t_n)$ of grade $0$ or $\omega$, where $s_i$ and $t_i\;
(1\le i\le n)$ are equivalent terms, then $\Tsr{\alpha}{\Gamma.}$ \\
\item If $\Ts{\beta}{\varrho}{\Gamma_i}$ and $\beta\lhd\alpha$ hold for every
premiss $\Gamma_i$ of an inference
$(\wedge),(\vee),\linebreak(\exists_1),(\forall_2)$ or $(Cut)$ with a cut formula
having grade $<\varrho$, and conclusion $\Gamma$, then
$\Tsr{\alpha}{\Gamma.}$\\
\item If $\Tsr{\al_0}{\Gamma,F(U)}$ holds for some $\al_0\lhd\al$
and a non-arithmetic formula $F(U)$ (i. e., $gr(F(U))\geq \om)$, then
$\Tsr{\al}{\Gamma,\exists X F(X)}.$
\item $(\omega$-rule$)$. If $\Ts{\beta}{\varrho}{\Gamma,A(\mq)}$ is true for
every
$m<\omega$, $\forall x A(x)\in \Gamma$, and $\beta\lhd\alpha$, then
$\Ts{\alpha}{\varrho}{\Gamma}.$\\
\item ($\Omega$-rule). Let $f$ be a fundamental function satisfying
\begin{enumerate}
\item[(a)] 
$f(\Om)\unlhd\alpha,$
\item[(b)]
$\Tsr{f(0)}{\Gamma,\forall X
F(X)}$, where $\forall X F(X)\in \Pi^1_1$, and
\item[(c)] $\Ts{\beta}{0}{\Xi,\forall X F(X)}$ implies
$\Tsr{f(\beta)}{\Xi,\Gamma}$ for every set of weak formulas $\Xi$ and
$\beta<\Om$.
\end{enumerate}
Then $\Tsr{\alpha}{\Gamma}$ holds.
\end{enumerate}}
\end{definition}
\begin{remark}{\em
The derivability relation $\Tsr{\alpha}{\Gamma}$ is from \cite{rathjen-weiermann93} and is modelled upon the
relation $\provx{PB^*}{\alpha}{n}{F} $ of \cite{bu},
the main difference being the sequent calculus setting instead of $P$--
and $N$--forms and a different assignment of cut--degrees. The allowance for
transfinite cut--degrees will enable us to deal with arithmetical comprehension.}
\end{remark}
\begin{remark}\label{Bedenken}{\em If one ruminates on the definition of
the derivability predicate $\Tsr{\alpha}{\Xi}$ the question arises whether it is actually a proper inductive definition. The critical point is obviously the condition (c) of the $\Omega$-rule. Note that $\Ts{\beta}{0}{\Xi,\forall X F(X)}$ occurs negatively
 in clause (c). However, since $\beta<\Omega$, the pertaining derivation does not contain any applications of the $\Omega$-rule. Thus
the definition of $\Tsr{\alpha}{\Xi}$ proceeds via an iterated inductive definition. First one defines a derivability
predicate without involvement of the $\Omega$-rule via an ordinary inductive definition, and in a second step
defines $\Tsr{\alpha}{\Gamma}$ inductively referring to the first derivability predicate in the $\Omega$-rule.

It will actually be a non trivial issue how to handle such inductive definitions in a weak background theory.
}\end{remark}
\begin{lemma}
\begin{enumerate}\label{6.4}
\item$\Ts{\alpha}{\delta}{\Gamma}\;\&\;\Gamma\subseteq\Delta\;\&\;
\alpha\unlhd\beta\;\&\;\delta\le\varrho
\;\Longrightarrow\;
\Tsr{\beta}{\Delta},$
\item$\Tsr{\alpha}{\Gamma,A\wedge B} \;\Longrightarrow\;
\Tsr{\alpha}{\Gamma,A}\;\&\;\Tsr{\alpha}{\Gamma,B,}$
\item$\Tsr{\alpha}{\Gamma,A\vee B} \;\Longrightarrow\;
\Tsr{\alpha}{\Gamma,A,B}$
\item$\Tsr{\alpha}{\Gamma,F(t)} \;\Longrightarrow\;
\Tsr{\alpha}{\Gamma,F(s)}$ if $t$ and $s$ are equivalent,\\
\item$\Tsr{\alpha}{\Gamma, \forall x F(x)} \;\Longrightarrow\;
\Tsr{\alpha}{\Gamma,F(s)}$ for every term $s$.\\
\item If $\Tsr{\alpha}{\Gamma,\forall X G(X)}$ and
$gr(G(U))\geq\omega$, then
$\Tsr{\alpha}{\Gamma,G(U)}$.
\end{enumerate}
\end{lemma}
\prf Proceed by induction on $\alpha$. These can be carried out
straightforwardly. (5) requires (4). As to (6), observe
that
$\forall XG(X) $ cannot be the main formula of an axiom.\hfill$\Box$
\begin{lemma}\label{6.5}
$\Ts{2\cdot \alpha}{0}{\Gamma,A(s_1,\ldots,s_k),\neg A(t_1,\ldots,t_k)}$
if $\alpha\geq gr(A(s_1,\ldots,s_k))$ and $s_i$ and $t_i$ are equivalent
terms.
\end{lemma}
\prf Proceed by induction on $gr(A(s_1,\ldots,s_k))$. Crucially note that if $gr(A(s_1,\ldots,s_k))=\omega$ then $\Gamma,A(s_1,\ldots,s_k),\neg A(t_1,\ldots,t_k)$ is an axiom according to Definition \ref{8.3} clause (4).
\qed

\begin{lemma}\label{6.6}
\begin{enumerate}
\item$\Ts{2m}{0}{\neg(0\in U),(\exists x)[x\in U\wedge \neg (x'\in
U)],\mq\in U},$
\item$\Ts{\omega+5}{0}{\forall X [0\in X\wedge \forall x(x\in
X\rightarrow x'\in X)\rightarrow \forall x (x\in X)].}$
\end{enumerate}
\end{lemma}
{\em Proof.} For (1) use induction on $m$.
(2) is an immediate consequence of (1) using
Lemma
\ref{6.4}
(1), the
$\omega$-rule, $(\vee)$, and $(\forall_2)$.
\begin{definition}{\em
For formulas $F(U)$ and $A(a)$, $F(A)$ denotes the result of replacing
each occurrence of the form $e\in U$ in $F(U)$ by $A(e)$. The
expression $F(A)$ is a formula if the bound variables in $A(a)$ are
chosen in an appropriate way, in particular, if $F(U)$ and $A(a)$ have
no bound variables\ in common.}
\end{definition}
\begin{lemma}\label{6.8}
 Suppose $\alpha<\Omega$ and let $\Delta(U)=\{F_1(U),\ldots,F_k(U)\}$ be a set of weak formulas such
that $U$ doesn't occur in $\forall X F_i(X)\;(1\le i\le k)$.
For an arbitrary formula $A(a)$ we then have:
$$\Ts{\alpha}{0}{\Delta(U)}\;\Longrightarrow\;\Ts{\Om+\alpha}{0}{\Delta(A)}.$$
\end{lemma}
\prf Proceed by induction on $\alpha$. Suppose $\Delta(U)$ is an axiom.
Then either $\Delta(A)$ is an axiom too, or $\Ts{\omega+\omega}{0}{\Delta(A)}$
can be
obtained through use of Lemma \ref{6.5}. Therefore
$\Ts{\Om+\alpha}{0}{\Delta(A)}$ by Lemma \ref{6.4} (1).
If $\Ts{\alpha}{0}{\Delta(U)}$ is the result of an inference, then this
inference must be different from $(\exists_2)$, $(Cut)$, and the $(\Om-rule)$ since $\Delta(U)$ consists of weak formulas, the derivation is cut-free and $\alpha<\Omega$. For the remaining possible inference rules
the assertion follows easily from the induction hypothesis. \qed
\begin{lemma}\label{stern}
Let $\Gamma,\forall X F(X)$ be a set of weak formulas.
If $\Ts{\al}{0}{\Gamma,\forall X F(X)}$
and $\al<\Om$, then $\Ts{\al}{0}{\Gamma, F(U)}$.
\end{lemma}
\prf Use induction on $\al$. Note that $\forall XF(X)$
cannot be a principal formula of an axiom, since $\exists X\neg F(X)$
does not surface in such a derivation. Also, due to $\al<\Om$, the derivation
doesn't involve instances of the $\Om$-rule. Therefore the proof is straightforward.
\hfill $\Box$\\[0.2cm]
The role of the $\Om$-rule in our calculus $T^*_{\!_Q}$ is enshrined in the next
lemma.
\begin{lemma}\label{6.9}
$\Ts{\Om\cdot 2}{0}{\exists X F(X),\neg F(A)} $ for every arithmetic
formula $F(U)$ and arbitrary formula $A(a)$.
\end{lemma}
\prf \setcounter{equation}{0}
Let $f(\alpha):=\Om+\alpha$ with $dom(f):=\{\alpha\in
\TP\!:\alpha\leq\Om\}.$ Then
\begin{equation}\label{6.1}
\Ts{f(0)}{0}{\forall X\neg F(X),\exists X F(X),\neg F(A)}
\end{equation}
according to Lemma \ref{6.5}.
For $\alpha<\Om$ and every set of weak formulas $\Theta$, we have by
Lemmata \ref{6.8} and \ref{stern},
$$\Ts{\alpha}{0}{\Theta,\forall X\neg F(X)}\;\Longrightarrow\;
\Ts{f(\alpha)}{0}{\Theta,\neg F(A).}$$
Therefore, by Lemma \ref{6.4} (1),
\begin{equation}\label{6.2}
\Ts{\alpha}{0}{\Theta,\forall X \neg F(X)}\;\Longrightarrow\;
\Ts{f(\alpha)}{0}{\Theta,\exists X F(X), \neg F(A).}
\end{equation}
The assertion now follows from (\ref{6.1}) and (\ref{6.2}) by the
$\Om$-rule.\hfill$\Box$

\begin{cor}\label{aca}
$\Ts{\Om\cdot 2+1}{\omega}{\exists X \,\forall y\,(y\in X\gdw B(y)) } $ for every arithmetic
formula $B(a)$.
\end{cor}
\prf Owing to Lemma \ref{6.9} we have
\begin{eqnarray}\label{aca1}&&\Ts{\Om\cdot 2}{0}{\exists X \,\forall y\,(y\in X\gdw B(y)),\,
 \neg \forall y\,(B(y)\gdw B(y))} .\end{eqnarray}
 As Lemma \ref{6.5} yields $\Ts{k}{0}{\forall y\,(B(y)\gdw B(y))}$ for some
 $k<\omega$, cutting with (\ref{aca1}) yields
  $\Ts{\Om\cdot 2+1}{\omega}{\exists X \,\forall y\,(y\in X\gdw B(x)) } $.
  \qed

\begin{cor}\label{6.9a} For every arithmetic relation $\prec$ (parameters allowed) and
arbitrary formula $A(a)$ we have
$\Ts{\Om\cdot 2+\omega}{0}{\forall \vec X\,\forall \vec x (\WF(\prec) \to \TI(\prec, A))} $
where the quantifiers $\forall \vec X\,\forall \vec x$ bind all free variables in  $\WF(\prec) \to \TI(\prec, A)$.
\end{cor}
\prf By Lemma \ref{6.9} we have $\Ts{\Om\cdot 2}{0}{\neg (\WF(\prec))',(\TI(\prec,A))'} $
where $'$ denotes any assignment of free numerical variables to numerals.
Hence $$\Ts{\Om\cdot 2+2}{0}{(\WF(\prec)\to\TI(\prec,A))'}$$ by two applications of $(\vee)$. Applying the $\omega$-rule
 the right number of times followed by the right number of $(\forall_2)$ inferences, one arrives at the desired conclusion.
 \qed

\subsection{The reduction procedure for $T^*_{\!_Q}$}
Below we follow \cite{rathjen-weiermann93} section 7.
\begin{lemma}\label{t7.1}
Let $C$ be a formula of grade $\varrho$. Suppose $C$ is a prime formula
or of either form $\exists X H(X),\;\exists x G(x)$ or $A\vee B$. Let
$\alpha=\omega^{\alpha_1}+\cdots+\omega^{\alpha_k}$ with $\delta\le
\omega^{\alpha_k}\le\ldots\le \omega^{\alpha_1}.$ Then we have
$\Tsr{\alpha}{\Delta,\neg
C}\;\&\;\Tsr{\delta}{\Gamma,C}\;\Longrightarrow
\Tsr{\alpha+\delta}{\Delta,\Gamma}.$
\end{lemma}
\prf We proceed by induction on $\delta$.\\
\setcounter{equation}{0}
1. Let $\Gamma,C$ be an axiom. Then there are three cases to consider.\\
1.1. $\Gamma$ is an axiom. Then so is $\Delta,\Gamma$. Hence
$\Tsr{\alpha+\delta}{\Delta,\Gamma}.$\\
1.2. $C$ is a true constant prime formula or negated prime formula. A
straightforward induction on $\alpha$ then yields
$\Tsr{\alpha}{\Delta}$, and thus $\Tsr{\alpha+\delta}{\Delta,\Gamma}$ by
\ref{6.4} (1).\\
1.3. $C\equiv A(s_1,\ldots,s_n)$ and $\Gamma$ contains a formula $\neg
A(t_1,\ldots,t_n)$ where $s_i$ and $t_i$ are equivalent terms. From
$\Tsr{\alpha}{\Delta,\neg A(s_1,\ldots,s_n)}$ one receives\\
$\Tsr{\alpha}{\Delta,\neg A(t_1,\ldots,t_n)}$ by use of Lemma \ref{6.4}
 (4).
Thence $\Tsr{\alpha+\delta}{\Delta,\Gamma}$ follows by use of Lemma
\ref{6.4}
 (1), since
$\neg A(t_1,\ldots,t_n)\in \Gamma$.\\
2. Suppose $C\equiv A\vee B$ and $\Tsr{\delta_0}{\Gamma,C,A_0}$ with
$\A_0\in\{A,B\}$ and $\delta_0\lhd \delta$. Inductively we get
\begin{equation}\label{7.1}
\Tsr{\alpha+\delta_0}{\Delta,\Gamma,A_0}.
\end{equation}
Next use Lemma \ref{6.4} (2) on
$\Tsr{\al}{\Delta,\neg A\wedge\neg B}$ to obtain
\begin{equation}\label{7.2}
\Tsr{\alpha+\delta_0}{\Delta,\Gamma,\neg A_0}.
\end{equation}
Whence use a cut on (\ref{7.1}) and (\ref{7.2}) to get the assertion.\\
3. Suppose $C\equiv\exists xG(x)$ and $\Tsr{\delta_0}{\Gamma,C,G(t)} $
with $\delta_0\lhd\delta$. Inductively we get
\begin{equation}\label{7.3}
\Tsr{\alpha+\delta_0}{\Delta,\Gamma,G(t)}.
\end{equation}
By Lemma \ref{6.4} 1), (5), we also get
\begin{equation}\label{7.4}
\Tsr{\alpha+\delta_0}{\Delta,\Gamma,\neg G(t)};
\end{equation}
thus (\ref{7.3}) and (\ref{7.4}) yield
$\Tsr{\alpha+\delta}{\Delta,\Gamma}$ by $(Cut)$.\\
4. Suppose the last inference was $(\exists_2)$ with p. f. $C$. Then
$C\equiv\exists X H(X)$ and $\Tsr{\delta_0}{\Gamma,C,H(U)} $
for some $\delta_0\lhd\delta$ and $gr(H(U))\geq\om$. Inductively we get
\begin{equation}\label{7.5}
\Tsr{\alpha+\delta_0}{\Delta,\Gamma,H(U).}
\end{equation}
By Lemma \ref{6.4}  (1), (6) we also get
\begin{equation}\label{7.6}
\Tsr{\alpha+\delta_0}{\Delta,\Gamma,\neg H(U).}
\end{equation}
From (\ref{7.5}) and (\ref{7.6}) we obtain
$$\Tsr{\alpha+\delta}{\Delta,\Gamma.}$$
5. Let $\Tsr{\delta}{\Gamma,C} $ be derived by the $\Omega$-rule with
fundamental function $f$. Then the assertion follows from the I. H. by
the $\Omega$-rule using the fundamental function ${\alpha}+f$.\\
6. In the remaining cases the assertion follows from the I. H. used on
the premises and by reapplying the same inference.\hfill$\Box$
\begin{lemma}\label{c7.2}
$\Ts{\alpha}{\eta+1}{\Gamma}\;\Longrightarrow\;\Ts{\omega^{\alpha}}{\eta}{\Gamma.}$

\end{lemma}
\prf We proceed by induction on $\alpha$. We only treat the crucial case
when $\Ts{\alpha_0}{\eta+1}{\Gamma,D}$ and
$\Ts{\alpha_0}{\eta+1}{\Gamma,\neg D}$, where $\alpha_0\lhd \alpha$,
and $gr(D)=\eta$.
Inductively this becomes
$\Ts{\omega^{\alpha_0}}{\eta}{\Gamma,D}$ and
$\Ts{\omega^{\alpha_0}}{\eta}{\Gamma,\neg D.}$
Since $D$ or $\neg D$ must be one of the forms exhibited in Lemma
\ref{t7.1}, we
obtain $\Ts{\omega^{\alpha_0}+\omega^{\alpha_0}}{\eta}{\Gamma}$ by Lemma
\ref{t7.1}.
 As $\omega^{\alpha_0}+\omega^{\alpha_0}\lhd \omega^{\alpha}$, we can use
Lemma \ref{6.4} {\em 1.}) to get the assertion.
\begin{theorem}[{\em Collapsing Theorem}]\label{kollaps}
Let $\Gamma$ be a set of weak formulas.
 We
have
$$\Ts{\alpha}{\omega}{\Gamma}\;\Longrightarrow\;\Ts{\vt\alpha}{0}{\Gamma.}$$
\end{theorem}
\prf We proceed by induction on  $\alpha$.\setcounter{equation}{0}
Observe that for $\be<\de<\Om$, we always have $\be\lhd\de.$\\
1. If $\Gamma$ is an axiom, then the assertion is trivial.\\
2. Let $\Ts{\alpha}{\omega}{\Gamma}$ be the result of an inference
other than $(Cut)$ and $\Omega$-rule. Then we have
$\Ts{\alpha_0}{\omega}{\Gamma_i}$ with $\alpha_0\lhd
\alpha$ and $\Gamma_i$ being the $i$-th premiss of that
inference.
$\al_0\lhd\al$  implies $\vt\al_0\lhd\vt\al$. Therefore
$\Ts{\vt\al_0}{0}{\Gamma_0}$ by the I. H., hence
$\Ts{\vt\alpha}{0}{\Gamma}$ by reapplying the same inference.\\
3. Suppose
$\Ts{\alpha}{\omega}{\Gamma}$ results by the $\Omega$-rule with respect
to a $\Pi^1_1$-formula $\forall X F(X)$ and a fundamental function $f$.
Then  $f(\Om)\unlhd \alpha$ and
\begin{equation}\label{731}
\Ts{f(0)}{\omega}{\Gamma,\forall X F(X),}
\end{equation}
and, for every set of weak formulas $\Xi$ and $\beta<\Om$,
\begin{equation}\label{732}
\Ts{\beta}{0}{\Xi,\forall X F(X)}\;\Longrightarrow\;
\Ts{f(\beta)}{\omega}{\Xi,\Gamma.}
\end{equation} The I. H. used on (\ref{731}) supplies us with
$\Ts{\vt (f(0))}{0}{\Gamma,\forall XF(X)}$. Hence with $\Xi=\Gamma$ we
get
\begin{equation}
\label{733}
\Ts{f(\vt(f(0)))}{\omega}{\Gamma}
\end{equation}
from (\ref{732}).
Now Lemma \ref{f5.4}  ensures that $f(\be)\lhd
f(\Om)$, where $\be=\vt(f(0))$.

So using the I. H. on (\ref{733}),
we obtain
\begin{equation}\label{734}
\Ts{\vt(f(\be))}{0}{\Gamma},
\end{equation}
thus $\Ts{\vt\al}{0}{\Gamma}$ as $f(\be)\lhd\al$.\\
4. Suppose
$\Ts{\alpha_0}{\omega}{\Gamma,A}$ and
$\Ts{\alpha_0}{\omega}{\Gamma,\neg A}$,
where $\alpha_0\lhd\alpha$ and $gr(A)<\omega$. Inductively we then get
$\Ts{\vt\alpha_0}{0}{\Gamma,A}$ and
$\Ts{\vt\alpha_0}{0}{\Gamma,\neg A.}$
Let $gr(A)=n-1$.
Then (Cut) yields
\begin{equation}
\Ts{\be_1}{n}{\Gamma}
\end{equation}
with $\be_1=(\vt\al_0)+1$.
Applying Lemmma \ref{c7.2}, we get
$\Ts{\om^{\beta_1}}{n-1}{\Gamma}$, and by
repeating this process we arrive at
$$\Ts{\beta_n}{0}{\Gamma},$$
where $\beta_{k+1}:=\omega^{\beta_k}\;(1\le k<n)$.
Since $\vt\al_0<\vt\al$, we have $\be_n<\vt\al;$
thus
$\Ts{\vt\alpha}{0}{\Gamma.}$ \qed

\subsection{Embedding $\MQ$  into  $T^*_{\!_Q}$.}
Assuming that $\MQ$ is well-founded tree, the objective of this section is to embed $\MQ$ into $T^*_{\!_Q}$, so as to
obtain a contradiction.
Let $\fXO$ be the Kleene-Brouwer ordering of $\MQ$.
We write $\provx{\MQ}{\tau}{}{\Gamma}$ if $\Gamma$ is the sequent attached to the node $\tau$ in $\MQ$.

\begin{theorem}\label{einbett}
$\provx{\MQ}{\tau}{}{\Xi} \Rightarrow \exists k<\omega\, \Ts{\EE_{\tau}+k}{\omega}{\Xi}$.

\end{theorem}
\prf We proceed by induction on $\tau$, i.e., the  Kleene-Brouwer ordering of $\MQ$.

Suppose $\tau$ is an end-node of $\MQ$. Then $\Xi$ must be axiomatic and therefore is an axiom of $T^*_{\!_Q}$,
and hence $\Ts{\EE_{\tau}}{\omega}{\Xi}$.

Now assume that $\tau$ is not an end-node  of $\MQ$. Then $\Xi$ is not axiomatic.

If $\Xi$ is not reducible, then there is a node $\tau_0$ immediately above $\tau$ in $\MQ$ such that
$\provx{\MQ}{\tau_0}{}{\Xi, \neg \bar Q(i),\neg A_i}$ for some   $i$. Inductively we have
$$\Ts{\EE_{\tau_0}+k_0}{\omega}{\Xi,\neg \bar Q(i),\neg A_i}$$ for some $k_0<\omega$. We also have $\Ts{0}{0}{\bar Q(i)}$ and, using Corollary \ref{aca} (if $i=0$) and Corollary \ref{6.9a} (if $i>0$),
$\Ts{\Omega\cdot 2+\omega}{\omega}{A_i}$. Thus, noting that $\Omega\cdot 2+\omega\lhd \EE_{\tau_0}+k_0$, and by employing two cuts we arrive at
$$\Ts{\EE_{\tau_0}+k_0+2}{\omega+n}{\Xi}$$ for some $n<\omega$.
By Lemma \ref{c7.2} we get $\Ts{\omega_n(\EE_{\tau_0}+k_0+2)}{\omega}{\Xi}$, and hence
$\Ts{\EE_{\tau}}{\omega}{\Xi}$ since $\omega_n(\EE_{\tau_0}+k_0+2)\lhd \EE_{\tau}$.

Now suppose that $\Xi$ is reducible.
$\Xi$ will be of the form
$$\Xi', \thinspace E,\thinspace \Xi''$$
where $E$ is not a literal and $\Xi'$ contains only literals.

First assume $E$ to be of the form $\forall x\,F(x)$.
Then, for each $m$, there is a node $\tau_m$ immediately above $\tau$ in $\MQ$ such that
$$\provx{\MQ}{\tau_n}{}{\Xi',F(\bar{m}),\Xi'', \neg \bar Q(i),\neg A_i}$$
for some $i$. Inductively we have
$$\Ts{\EE_{\tau_m}+k_m}{\omega}{\Xi',F(\bar{m}),\Xi'',\neg \bar Q(i),\neg A_i}$$ for all $m$, where $k_m<\omega$.
We also have $\Ts{0}{0}{\bar Q(i)}$ and, using Lemma \ref{6.9a},
$\Ts{\Omega\cdot 2+\omega}{0}{A_i}$. Thus, noting that $\Omega\cdot 2+\omega\lhd \EE_{\tau_m}+k_m$, and by employing two cuts
there is an $n$ such that
$$\Ts{\EE_{\tau_m}+k_m+2}{\omega+n}{\Xi',F(\bar{m}),\Xi''}$$
holds for all $m$. By Lemma \ref{c7.2} we get
$$\Ts{\omega_n(\EE_{\tau_m}+k_m+2)}{\omega}{\Xi',F(\bar{m}),\Xi''}$$
 for all $m$. Whence
 $$\Ts{\EE_{\tau}}{\omega}{\Xi',F(\bar{m}),\Xi''}$$ since $\omega_n(\EE_{\tau_m}+k_m+2)\lhd \EE_{\tau}$.
 A final application of the
$\omega$-rule yields
$$\Ts{\EE_{\tau}+1}{\omega}{\Xi',\forall x\,F(x),F(\bar{m}),\Xi''}$$
i.e.,  $\Ts{\EE_{\tau}+1}{\omega}{\Xi}$.

If $E$ is a redex of another type but not of the form $\exists X B(X)$ with
 $B(U)$ arithmetic, then one proceeds in a similar way as in the previous case.

Now assume $E$ to be of the form $\exists X\,B(X)$ with $B(U)$ arithmetic.
Then there is a node $\tau_0$ immediately above $\tau$ in $\MQ$ such that
$$\provx{\MQ}{\tau_0}{}{\Xi',B(U),\Xi'', \neg \bar Q(i),\neg A_i}$$
for some $i$ and set variable $U$. Inductively we have
$$\Ts{\EE_{\tau_0}+k_0}{\omega}{\Xi',B(U),\Xi'',\neg \bar Q(i),\neg A_i}$$ for some $k_0<\omega$.
We also have $\Ts{0}{0}{\bar Q(i)}$ and, using Lemma \ref{6.9a},
$\Ts{\Omega\cdot 2+\omega}{0}{A_i}$. Thus, noting that $\Omega\cdot 2+\omega\lhd \EE_{\tau_0}+k_0$, and by employing two cuts
there is an $n$ such that
$$\Ts{\EE_{\tau_0}+k_0+2}{\omega+n}{\Xi',B(U),\Xi''.}$$
By Lemma \ref{c7.2} we get
\begin{eqnarray}\label{10.1.1} &&\Ts{\omega_n(\EE_{\tau_0}+k_0+2)}{\omega}{\Xi',B(U),\Xi''.}\end{eqnarray}
Lemma \ref{6.9} yields
\begin{eqnarray}\label{10.1.2} &&\Ts{\Om\cdot 2}{0}{\exists X B(X),\neg B(U).}\end{eqnarray}
Cutting $B(U)$ and $\neg B(U)$ out of (\ref{10.1.1}) and (\ref{10.1.2})
we arrive at
$$\Ts{\omega_n(\EE_{\tau_0}+k_0+2)+1}{\omega}{\Xi',\exists X B(X),\Xi''.}$$
Since $\omega_n(\EE_{\tau_0}+k_0+2)+1\lhd \EE_{\tau}$ we get
$\Ts{\EE_{\tau}}{\omega}{\Xi',\exists X B(X),\Xi''}$,
 i.e.,  $\Ts{\EE_{\tau}}{\omega}{\Xi}$.
  \qed

Below $\emptyset$ stands for the empty sequent and $\tau_0$ denotes the bottom node of $\MQ$ which is the maximum element
of the pertaining Kleene-Brouwer ordering.
\begin{cor}\label{end1} If $\MQ$ is well-founded, then $\Ts{\vt(\omega_n(\EE_{\tau_0}+m))}{0}{\emptyset}$ for some $n,m<\omega$.
\end{cor}
\prf We have $\provx{\MQ}{\tau_0}{}{ \neg \bar Q(0),\neg A_0}$.
Thus there is a $k<\omega$ such that $$\Ts{\EE_{\tau_0}+k}{\omega}{\neg \bar Q(0),\neg A_0}$$  holds by Theorem \ref{einbett}.
 We also have $\Ts{0}{0}{\bar Q(0)}$ and, using Corollary \ref{6.9a},
$\Ts{\Omega\cdot 2+\omega}{0}{A_0}$. Thus, noting that $\Omega\cdot 2+\omega\lhd \EE_{\tau_0}+k$, and by employing two cuts we arrive at
$$\Ts{\EE_{\tau_0}+k+2}{\omega+n}{\emptyset}$$ for some $n<\omega$.
Via Lemma \ref{c7.2} we deduce $\Ts{\omega_n(\EE_{\tau_0}+k+2)}{\omega}{\emptyset}$, so that by Theorem \ref{kollaps}
we conclude $\Ts{\vt(\omega_n(\EE_{\tau_0}+m))}{0}{\emptyset}$ with $m=k+2$. \qed

\begin{cor}\label{end2} $\MQ$ is not well-founded. \end{cor}
\prf If $\MQ$ were well-founded we would have \begin{eqnarray}\label{end21}&&\Ts{\vt(\omega_n(\EE_{\tau_0}+m))}{0}{\emptyset}\end{eqnarray} for some $n,m<\omega$
by Corollary \ref{end1}. But a straightforward induction on $\alpha<\Omega$ shows
that  $$\Ts{\alpha}{0}{\Gamma}\;\Rightarrow \;\Gamma\ne\emptyset,$$
yielding that (\ref{end21}) is impossible.\qed

It remains to show that the result of Corollary \ref{end2} is provable in $\ACA_0$ from
$$\forall {\mathfrak X}\,(\WO({\mathfrak X})\to \WO(\vt_{\mathfrak X}))\,.$$
Let $\Sr$ be the theory $\ACA_0$ plus the latter axiom.
The main issue is how to formalize the derivability predicate
$\Ts{\alpha}{\rho}{\Gamma}$ in the background theory $\Sr$. We elaborated earlier in Remark  \ref{Bedenken} that this seems to require
an iterated inductive definition, something apparently not available in $\Sr$.
However, all we need is a fixed point not a proper inductive definition, i.e., to capture the notion of derivability in $\Tst$
without the $\Omega$-rule it suffices to find a predicate $\cal D$ of $\alpha,\rho,\Gamma$ such that
\begin{itemize}\item[$(*)$] \label{D}
$\calD(\alpha,\rho,\Gamma)$ if and only if  $\alpha\in\vert \vartheta_{\mathfrak X}\vert$, $\rho\leq\omega+\omega$, $\Gamma$ is a sequent, and  either $\Gamma$ contains an axiom of $\Tst$
or $\Gamma$ is the conclusion
of an inference of $\Tst$ other than $(\Omega)$ with premisses $(\Gamma_i)_{i\in I}$ such that for every $i\in I$ there exists $\beta_i\lhd\alpha$ with $\calD(\beta_i,\rho,\Gamma_i)$, and if the inference is a cut it has rank $<\rho$.
\end{itemize}
$(*)$ can be viewed as a fixed-point axiom which together with transfinite induction for  $\vartheta_{\mathfrak X}$
  defines $\Tst$-derivability (without $(\Omega)$-rule) implicitly.

  How can we find a fixed point as described in $(*)$? As it turns out, it follows from \cite{H-Band} that $\Sr$ proves that
  every set is contained in a countable coded $\omega$-model of the theory $\ATR_0$. It is also known that
  $\ATR_0$ proves the $\Sigma^1_1$ axiom of choice, $\SAC$ (see \cite[Theorem V.8.3]{SOSA}).
  Moreover, in $\ACA_0+\SAC$ one can prove  for every $P$-positive arithmetical formula $A(u,P)$ that there is
  a $\Sigma^1_1$ formula $F(u)$ such that $\forall x[F(x)\leftrightarrow A(x, F)]$, where
  $A(x, F)$ arises from $A(x,P)$ by replacing every occurrence of the form $P(t)$ in the first formula by $F(t)$.
  This is known as the Second Recursion Theorem (see \cite[V.2.3]{ba}).
  Arguing in $\Sr$, we find a countable coded $\omega$ model $\mathfrak B$ with $\mathfrak X\in\mathfrak B$ such that
  $\mathfrak B$ is a model of  $\ATR$. As a result, there is a predicate $\mathcal D$ definable in $\mathfrak B$ that satisfies
  $(*)$. As a result, $\mathcal D$ is a set in $\Sr$.
  To obtain the full derivability relation  $\Ts{\alpha}{\rho}{\Gamma}$ we have to take the $\Omega$-rule into account.
  We do this by taking a countable coded $\omega$-model $\mathfrak C$ of $\ATR$ that contains both $\mathfrak X$ and $\mathcal D$.
  We then define an appropriate fixed point predicate $\mathcal D_{\Omega}$ using the clauses for defining $\Ts{\alpha}{\rho}{\Gamma}$
  and $\mathcal D$ for the negative occurrences in the $\Omega$-rule.

  The upshot is that we can formalize all of this in $\Sr$.

\begin{remark}{\em
When giving talks about the material of this article, the first author was
asked what the proof-theoretic ordinal of the theories that Theorem \ref{main} is concerned with might be.
He conjectures  that it is the ordinal
$$\vartheta(\varphi2(\Omega+1))$$ (or $\psi (\varphi2(\Omega+1))$ in the
representation system based on the $\psi$-function; see \cite[section 3]{rathjen-weiermann93}), i.e. the collapse of the first fixed point of
the epsilon function above $\Omega$.}
\end{remark}

\paragraph{Acknowledgement:}
The first author acknowledges
 support by the EPSRC of the UK through
 grants  EP/G029520/1 and EP/G058024/1.

The authors would also like to thank an anonymous referee for very helpful comments
and numerous suggestions. We also thank Anton Setzer for comments on a draft version
of this paper.

The results of this article were incorporated in the PhD thesis \cite{thesis} of the second author.

\end{document}